\documentclass[11pt, reqno]{amsart}

\usepackage{amsmath}
\usepackage{amsxtra}
\usepackage{amscd}
\usepackage{amsthm}
\usepackage{amsfonts}
\usepackage{amssymb}
\usepackage{eucal}

\newcommand{\nn}{\nonumber}
\newcommand{\bea}{\begin{eqnarray}}
\newcommand{\ena}{\end{eqnarray}}
\newcommand{\be}{\begin{eqnarray*}}
\newcommand{\en}{\end{eqnarray*}}

\newcommand{\bra}[1]{\langle #1 |}        
\newcommand{\ket}[1]{{| #1 \rangle}}      

\newcommand{\Z}{{\mathbb Z}}  
\newcommand{\C}{{\mathbb C}}

\newtheorem{thm}{Theorem}[section]
\newtheorem{prop}[thm]{Proposition}
\newtheorem{lem}[thm]{Lemma}
\newtheorem{cor}[thm]{Corollary}

\numberwithin{equation}{section}

\begin{document}


\pagestyle{myheadings}

\title[Fermionic formulas for $(k, 3)$-admissible configurations]
{Fermionic formulas for $(k, 3)$-admissible configurations}


\author{B. Feigin, M. Jimbo, T. Miwa, E. Mukhin and Y. Takeyama}
\address{BF: Landau institute for Theoretical Physics, Chernogolovka,
142432, Russia}\email{feigin@feigin.mccme.ru}  
\address{MJ: Graduate School of Mathematical Sciences, University of
Tokyo, Tokyo 153-8914, Japan}\email{jimbomic@ms.u-tokyo.ac.jp}
\address{TM: Department of Mathematics, Graduate School of Science, 
Kyoto University, Kyoto 606-8502 Japan}\email{tetsuji@kusm.kyoto-u.ac.jp}
\address{EM: Department of Mathematics, 
Indiana University-Purdue University-Indianapolis, 
402 N.Blackford St., LD 270, 
Indianapolis, IN 46202}\email{mukhin@math.iupui.edu} 
\address{YT: Department of Mathematics, Graduate School of Science, 
Kyoto University, Kyoto 606-8502 Japan}\email{takeyama@kusm.kyoto-u.ac.jp}


\date{\today}

\setcounter{footnote}{0}\renewcommand{\thefootnote}{\arabic{footnote}}


\begin{abstract} 
We obtain the fermionic formulas for the characters of 
$(k, r)$-admissible configurations in the case of $r=2$ and $r=3$. 
This combinatorial object appears as a label of a basis 
of certain subspace $W(\Lambda)$ of 
level-$k$ integrable highest weight module of 
$\widehat{\mathfrak{sl}}_{r}$. 
The dual space of $W(\Lambda)$ is embedded 
into the space of symmetric polynomials. 
We introduce a filtration on this space 
and determine the components of the associated graded space 
explicitly by using vertex operators. 
This implies a fermionic formula for the character of $W(\Lambda)$. 
\end{abstract}


\maketitle



\setcounter{footnote}{0}
\renewcommand{\thefootnote}{\arabic{footnote})}
\renewcommand{\arraystretch}{1.2}

\setcounter{section}{0}
\setcounter{equation}{0}

\section{Introduction} 

Let $\widehat{\mathfrak{sl}}_{r}$ be the affine Lie algebra 
$\mathfrak{sl}_{r} \otimes \C[t, t^{-1}] \oplus \C c \oplus \C d$ 
and $L(\Lambda)$ the integrable highest weight module for 
dominant integral weight $\Lambda$ of level $k$. 
We denote by $\mathfrak{a}$ the commutative Lie subalgebra of 
$\widehat{\mathfrak{sl}}_{r}$ generated by the elements 
\bea 
e_{21}[n], e_{31}[n], \ldots , e_{r1}[n], \quad n \in \Z. \nn
\ena 
Consider the $\mathfrak{a}$-submodule 
\bea 
W(\Lambda) := U(\mathfrak{a}) v_{\Lambda}, 
\label{principal-space} 
\ena 
where $v_{\Lambda} \in L(\Lambda)$ is 
the highest weight vector satisfying $e_{ij}[n]v_{\lambda}=0, (n>0)$. 
Our problem is to find some formulas for the character of $W(\Lambda)$. 
In \cite{P}, Primc constructed a basis of $W(\Lambda)$. 
His basis consists of vectors parametrized by the combinatorial object 
called $(k, r)$-admissible configurations. 
We can introduce some degrees on $(k, r)$-admissible configurations 
and define the character, which is equal to that of $W(\Lambda)$. 
In \cite{bosonic1} certain formulas, called `bosonic formulas', 
for the character of $(k, r)$-admissible configurations are obtained (see also \cite{bosonic2}).
Connections to Jack and Macdonald polynomials are discussed in
\cite{jack, macdonald, root-of-unity}.
In this paper we find another type of formulas in the cases of $r=2$ and $r=3$. 

We start from an algebra $E_{\Lambda}$ isomorphic to $W(\Lambda)$ as vector spaces. 
The algebra $E_{\Lambda}$ is constructed by generators 
$e_{21}[n], \ldots , e_{r1}[n], (n \le 0)$ with some relations. 
We will obtain the Gordon-type (or `fermionic') formulas for the character of $E_{\Lambda}$ 
by using vertex operators.  

Let ${\mathcal W}$ be a vector space with non-degenerate
quadratic form $\langle\cdot,\cdot\rangle$, and let 
$\Gamma$ be 
an integral lattice in $\mathcal W$, i.e.,
$\langle\gamma_1,\gamma_2\rangle\in{\bf Z}$ for any
$\gamma_1,\gamma_2\in\Gamma$. With such data we can associate
a lattice vertex operator algebra $\mathcal V_\Gamma$.
The algebra $\mathcal V_\Gamma$ is generated by vertex operators $V(\gamma,z)$
$(\gamma\in\Gamma, z\in{\bf C})$. Let us take a set $\{p_1,\ldots,p_n\}$
of linearly independent vectors from $\mathcal W$, and consider the subalgebra
$\mathcal C$ generated by the vertex operators $a_1(z),\ldots,a_n(z)$
where $a_i(z)=V(p_i,z)$. The
operators $a_i(z)$ satisfy quadratic relations.
It can be easily formulated in the case when $\langle p_i, p_j\rangle\geq0$
for all $i,j$. In this case, we have
\begin{eqnarray}
&&[a_\alpha(z),a_\alpha(w)]_\pm=0,
\hbox{ where }+(\hbox{resp., }-)\hbox{ if }
\langle p_\alpha, p_\alpha\rangle
\hbox{ is odd (resp., even)},\label{C1}\\[0pt]
&&[a_\alpha(z),a_\beta(w)]=0,\hbox{ if }\alpha\neq\beta,\label{C2}\\
&&a_\alpha(z)\partial_z^la_\beta(z)=0\hbox{ for }
l<\langle p_\alpha, p_\beta\rangle.\label{C3}
\end{eqnarray}
If $\langle p_\alpha, p_\beta\rangle<0$, then the 
relations are also quadratic, but $a_\alpha(z)$ and $a_\beta(w)$ are 
not commutative. 

Let us continue the discussion under the condition
$\langle p_i, p_j\rangle\geq0$. 
It is important that the relations \eqref{C1}--\eqref{C3} are the set of
defining relations. This fact actually is equivalent to the following
statement about representations of $\mathcal C$. Let $a_\alpha[i]$
be components of $a_\alpha(z)$, i.e., 
$a_\alpha(z)=\sum_{i\in\Z} a_\alpha[i]z^i$.
Choose the irreducible representation of ${\mathcal V}_\Gamma$ with
the vacuum vector $v$ satisfying $a_\alpha[i]v=0$ for $i\leq 0$.
Consider the space $W={\mathcal C}v$.
Let $\theta:W\rightarrow{\bf C}$ be a linear functional. 
Define the function
\[
\Psi^{\alpha_1,\ldots,\alpha_m}_\theta(z_1,\ldots,z_m)
=\langle\theta,a_{\alpha_1}(z_1)\cdots a_{\alpha_m}(z_m)v\rangle.
\]
It has a form
\bea
\Psi^{\alpha_1,\ldots,\alpha_m}_\theta(z_1,\ldots,z_m)
=F(z_1,\ldots,z_m)\prod_iz_i\prod_{i<j}(z_i-z_j)^
{\langle p_{\alpha_i}, p_{\alpha_j}\rangle}, 
\label{C4}
\ena	
where $F$ is a polynomial which is
symmetric with respect to the transposition of $z_i$ with $z_j$
if $\alpha_i=\alpha_j$.

Let $S$ be the space of functions of the form (\ref{C4}).
More precisely, we have a direct sum
$S=\oplus S_{\alpha_1,\ldots,\alpha_m}$,
where the set of indices $(\alpha_1\ldots,\alpha_m)$ is defined up to
permutations. The statement is that the map $W^*\rightarrow S$
is an isomorphism. This fact is equivalent to the relations
\eqref{C1}--\eqref{C3}, and also gives a possibility of writing down
the character of the space $W$. The space $W$ is naturally graded
by $\mathop{\rm deg}a_\alpha[i]=i$ as well as the function space $S$
by $\mathop{\rm deg}z_i=1$, and we have the equality of the corresponding 
characters:
$\mathop{\rm ch}W=\sum_{(\alpha_1,\ldots,\alpha_m)}
                  \mathop{\rm ch}S_{\alpha_1,\ldots,\alpha_m}$. 
We have
\[
\mathop{\rm ch}S_{\alpha_1,\ldots,\alpha_m}
=\sum_{m_1,\cdots,m_n\ge 0}
\frac{q^{\sum_{j}m_{j}+\sum_{i<j}{\langle p_{\alpha_i}, p_{\alpha_j}\rangle}m_{i}m_{j}}}
{(q)_{m_1}\cdots(q)_{m_n}},
\]
where $m_j$ is the number of $i$ such that $\alpha_i=j$. Using this
formula, we get a Gordon-type formula for the character of the space $W$.

Our strategy is to compare the more complicated algebras with the algebras
like $\mathcal C$. Let us consider the simplest example.

In the algebra $E={\bf C}[e[0], e[-1], e[-2], \ldots]$ there are a sequence of ideals
$E\supset J_1\supset J_2\cdots$. Here $J_s$ is the ideal generated by the
components $e^{(s)}[i]$ of the current
$e(z)^s=(\sum e[i]z^i)^s=\sum e^{(s)}[i]z^i$.
We want to study the quotient $E/J_{k+1}=E_k$. Filter $E_k$ by ideals
\[
E_k\supset J_k\supset J_k^2\supset\cdots,
\]
and construct the corresponding associated graded algebra
\[
E^{(1)}_k=E_k/J_k\oplus J_k/J_k^2\oplus J_k^2/J_k^3\oplus\cdots.
\]
We denote by the same symbol $J_k$ the image of $J_k\subset E$ in $E_k$.
Note that $E_k/J_k\simeq E_{k-1}$ and the algebra $E^{(1)}_k$ is generated
by the components of the currents $e(z)\in E_k/J_k$ and
$e^{(k)}(z)\in J_k/J_k^2$. The current $e^{(k)}(z)$ corresponds to the
current $e(z)^k$ from $E$. In the algebra $E_k/J_k\simeq E_{k-1}$
we have ideals 
$E_{k-1}\supset J_1\supset J_2\supset\cdots \supset J_{k-1}$.
Let $J^{(1)}_s$ be an ideal in $E^{(1)}_k$ generated
by $J_s\subset E_{k-1}$. In $E^{(1)}_k$ there are ideals
$J^{(1)}_1\supset J^{(1)}_2\supset\cdots\supset J^{(1)}_{k-1}$.
We can repeat such a construction and get the algebra
\[
E^{(2)}_k=E^{(1)}_k/J^{(1)}_{k-1}\oplus J^{(1)}_{k-1}/(J^{(1)}_{k-1})^2\oplus
\cdots.
\]
Obviously, the algebra $E^{(2)}_k$ is generated by the component of
the currents $e(z)$, $e^{(k)}(z)\in E^{(1)}_k/J^{(1)}_{k-1}$ and
$e^{(k-1)}(z)\in J^{(1)}_{k-1}/(J^{(1)}_{k-1})^2$. In
$E_{k-2}=E^{(1)}/J^{(1)}_{k-1}$ we have its own sequence of ideals
$E_{k-2}\supset J_1\supset J_2\supset\cdots\supset J_{k-2}$, and
we can repeat what we did before. As a result we get an algebra $E^{(2)}_k$
which is generated by $e(z),e^{(k-1)}(z),e^{(k)}(z)$. 
Then, filter $E^{(2)}_k$ again, and so on.  
In the end we construct an algebra $E^{(k)}_k$, which we denote by 
$B_k$. The algebra $B_k$ is generated by
$e(z)=e^{(1)}(z),e^{(2)}(z),\ldots,e^{(k)}(z)$. It has many gradings.
Surely, it inherits the $q$-grading, ${\rm deg}e_i=i$. It has also
${\bf Z}^k$-grading: each of the generators $e^{(\alpha)}(z)$
is homogeneous and has grading
$(0,\ldots,\buildrel{\alpha-{\rm th}}\over1,\ldots,0)$.
By a simple calculation, one can check that the generators $e^{(\alpha)}(z)$
satisfy the quadratic relations,
\bea
e^{(\alpha)}(z) \partial_z^l e^{(\beta)}(z)=0
\hbox{ for }l<2 \min(\alpha,\beta).
\label{QR}
\ena
Actually, these relations are defining relations for $B_k$. One way to prove
this statement is to compare $B_k$ with some algebra generated by vertex
operators. Now, we explain how to do it.

Consider an integrable representation of $\widehat{{\mathfrak{sl}}}_2$
of level $k$. It is known that in such a representation 
the current $e_{21}(z)$ 
satisfies the relation $e_{21}(z)^{k+1}=0$. 
Here $e_{21}$ is the nilpotent generator
of ${\mathfrak{sl}}_2$ and $e_{21}(z)$ is the corresponding current.

The explicit construction of such an 
$e_{21}(z)$ uses the so-called vertex operator realization. 
To do it consider the vector space ${\mathcal W}$ with a base
$ p_1,\ldots, p_k$ and a bilinear form
$\langle p_i,p_j\rangle=2\delta_{i,j}$. Let $a_i(z)=V(p_i,z)$ and
$b_i(z)=V(-p_i,z)$. Let $e_{21}(z)=a_1(z)+\cdots+a_k(z)$ and
$e_{12}(z)=b_1(z)+\cdots+b_k(z)$. It is well-known that such $e_{21}(z)$ and $e_{12}(z)$
generate $\widehat{{\mathfrak{sl}}}_2$ of level $k$.
The whole construction is nothing but the tensor product of $k$ copies
of the standard vertex oeprator realization of $\widehat{{\mathfrak{sl}}}_2$
of level $1$.

The representation of the corresponding vertex operator algebra after
restriction to $\widehat{{\mathfrak{sl}}}_2$ is a sum of integrable
representations of level $k$. Choose the vacuum vector $v$ in the
representation $\mathcal F$ of the vertex operator algebra which generates
the vacuum module for $\widehat{{\mathfrak{sl}}}_2$. (Our convention is such
that $e_{ij}[n]v=0$ for $n>0$.) There is a map $\varphi:E_k\rightarrow\mathcal F$
such that $P(e[0],e[-1],\ldots)\buildrel\varphi\over\mapsto P(e[0], e[-1], \ldots)v$.
We will prove that $\varphi$ is an embedding.

\def\e{\varepsilon}

Consider the family 
of maps $\varphi_\varepsilon:E_k\rightarrow\mathcal F$
where $\e\in{\bf C}, \e\not=0$ which send $e_{21}[i]$ to the $i$-th component of the
current $e_\e(z)=a_1(z)+\e a_2(z)+\cdots+\e^{k-1}a_k(z)$.
Let $\varphi_0$ be the limit of $\varphi_\e$ when $\e\rightarrow0$.
More precisely, we want to study the limit $W_0$
of the image of $\varphi_\e$ when $\e\rightarrow0$. 
First consider the limit of operators
\begin{eqnarray*}
&&\lim_{\e\rightarrow0}e_\e(z)=a_1(z),
\\
&&\lim_{\e\rightarrow0}\e^{-1}e_\e(z)^2=a_1(z)a_2(z),
\\
&& \qquad \ldots, 
\\
&&\lim_{\e\rightarrow0}\e^{1-s}e_\e(z)^s=a_1(z)\cdots a_s(z).
\end{eqnarray*}
Note that $\rho_s(z)=a_1(z)\cdots a_s(z)$  are vertex operators
$V(q_s,z)$ where $q_{s}$'s are vectors such that 
$\langle q_\alpha,q_\beta\rangle=2{\rm min}(\alpha,\beta)$.
It means that they satisfy the same quadratic relations as generators
$e_{21}^{(\alpha)}(z)$ in the algebra $B_k$. 
Looking more carefully at the limit 
$\e\rightarrow0$, it is possible to show that we have a surjection
$B_k\rightarrow{\rm lim}_{\e\rightarrow0}\varphi_\e(E_k)$. It means
that there is a family of algebras $U_\e$ such that $U_0\simeq B_k$, 
$U_\e\simeq E_k$ for $\e\not=0$, and
$\varphi_\e:U_\e\rightarrow\mathcal F$.
Therefore, we have a surjection $B_k\rightarrow W_0$, and in $W_0$
there is a subspace $\widetilde W_0={\mathcal C}v$ where $\mathcal C$
is the algebra generated by $\rho_s(z)$.

Comparing the characters of $B_k,\widetilde W_0$ and $W_0$ 
we get that
actually they are all isomorphic. As a corollary, we establish the Gordon-type
formula for the character of $E_k$.

There are many cases that can be studied in a similar manner.
We can replace 
$\widehat{{\mathfrak{sl}}}_2$ by $\widehat{\mathfrak g}$ for
any simply-laced semi-simple Lie algebra $\mathfrak g$.
Let $\mathfrak n$ be a 
maximal nilpotent subalgebra in $\mathfrak g$,
$L_k$ the vacuum representation of $\widehat{\mathfrak g}$ of level $k$,
$v$ the vacuum vector of $L_k$ and $W_0=U(\widehat{\mathfrak n})v\subset L_k$.

Following \cite{FK} we can realize $L_1$ as a representation of
some lattice vertex operator algebra. In this construction,
simple root generators $g_\alpha(z)\in\widehat{\mathfrak n}$
are just vertex operators
(up to some twisting, which is not essential in our argument).
Operators $g_\alpha(z)$ in representation of level $k$ can be represented
as a sum of vertex operators:
$g_\alpha(z)=g^{(1)}_\alpha(z)+g^{(2)}_\alpha(z)+\cdots+g^{(k)}_\alpha(z)$.
Now let us use the same $\e$-method.
Namely, introduce operators 
\[
g_{\alpha,\e}(z)=g^{(1)}_\alpha(z)+\e g^{(2)}_\alpha(z)+\cdots
+\e^{k-1}g^{(k)}_\alpha(z).
\]
Again, we consider the limit $\e\rightarrow0$, and repeating the process
in the $\widehat{{\mathfrak sl}}_2$ case we get the following result, 
which was formulated in \cite{FS}.
\[
\mathop{\rm ch}W_0=\sum_{m_1,\ldots,m_{kr}}
\frac{q^{\frac12\langle D{\bf m},{\bf m}\rangle}}
{(q)_{m_1}\cdots(q)_{m_{kr}}}.
\]
Here $r$ is the rank of $\mathfrak g$ and $D$ is the 
tensor product of two matrices
$C\otimes G$, where $C$ is the Cartan matrix of $\mathfrak g$ and
$G$ is the $k\times k$ Gordon matrix, i.e., $G=(G_{\alpha,\beta})$
where $G_{\alpha,\beta}={\rm min}(\alpha,\beta)$.

In a slightly different manner, the same method is used for the problem
which we will consider in this paper.

Let $L_k$ be the vacuum representation of $\widehat{{\mathfrak{sl}}}_3$ 
of level $k$. By $e_{ij}(z)$ we denote the standard basis of
$\widehat{{\mathfrak{sl}}}_3$ ($1\leq i,j\leq 3$).
Let $a(z)=e_{21}(z)$ and $b(z)=e_{31}(z)$. It is known that in $L_k$
the currents $a(z)$ and $b(z)$ satisfy the relations
\[
a(z)^\alpha b(z)^\beta=0\hbox{ if }\alpha+\beta=k+1.
\]
These are equivalent to the integrability of representation.
For $k=1$, $a(z)$ and $b(z)$ can be realized by vertex operators:
$a(z)=V(q_1,z)$ and $b(z)=V(q_2,z)$ where
$\langle q_1,q_1\rangle=\langle q_2,q_2\rangle=2$
and $\langle q_1,q_2\rangle=1$.
Again, for a bigger $k$, we consider
$a(z)=V(t_1,z)+\cdots+V(t_k,z)$ and $b(z)=V(s_1,z)+\cdots+V(s_k,z)$ where
$t_1,\ldots,t_k$ and $s_1,\ldots,s_k$ are vectors with the scalar products
$\langle t_i,t_i\rangle=\langle s_i,s_i\rangle=2$,
$\langle s_i,t_i\rangle=1$ and$
\langle t_i,t_j\rangle=\langle s_i,s_j\rangle
=\langle s_i,t_j\rangle=0$ for $i\not=j$.

Degeneration is given by formulas
\begin{eqnarray*}
a(\e,z)=V(t_1,z)+\e V(t_2,z)+\cdots+\e^{k-1}V(t_k,z),\\
b(\e,z)=\e^{k-1}V(t_1,z)+\e^{k-2}V(t_2,z)+\cdots+V(t_k,z).
\end{eqnarray*}
 The Gordon-type formula for the character of $W_0=\mathcal{C}v$
where $\mathcal{C}$ is generated by $a(z),b(z)$ and $v$ is the vacuum vector
in $L_k$, can be found in Theorem \ref{fermionic-formula:2k}. 
Again we have vertex operators which represent the currents 
$a(z)^\alpha$ and $b(z)^\beta$, $\alpha,\beta\le k$. 
Now $a(z)^\alpha a(w)^{\beta}\sim (z-w)^{2\min(\alpha,\beta)}$.
By this we mean that in the representation an arbitrary matrix element
$\langle\theta^{\vee},a(z)^\alpha a(w)^\beta\theta\rangle$ 
has the form $(z-w)^{2\min(\alpha,\beta)}f(z,w)$ where 
$f(z,w)$ is a Laurent polynomial. 
The currents $b(z)^\beta$ have the same properties and 
$a(z)^\alpha b(w)^\beta\sim (z-w)^{(\alpha+\beta-k)_+}$,
with $(m)_+=\max(m,0)$. 

In our paper we use the same $\varepsilon$-method to study 
$\widehat{\mathfrak{sl}}_3$ representations a little differently. 
Let us combine the currents $a(z)$ and $b(z)$  into a single one as
$e(z)=a(z^2)+zb(z^2)$. 
The relations $a^\alpha(z)b^\beta(z)=0$ ($\alpha+\beta=k+1$)
can be written in terms of $e(z)$ as 
$e^\alpha(z)e^\beta(-z)=0$  ($\alpha+\beta=k+1$).
The algebra
$\widehat{\mathfrak{sl}}_3$ has a vertex operator realization 
where the current $e(z)$ is a sum of vertex operators, 
and all previous techniques can be used. 
For $k=1$, $e(z)$ satisfies the relations 
$e(z)^2=0$ and $e(z)e(-z)=0$. The matrix elements
$\langle\theta^{\vee},e(z_1)e(z_2)\theta\rangle$ 
have the form $(z_1-z_2)^2(z_1+z_2)f(z_1,z_2)$ where $f(z_1,z_2)$ 
is a Laurent polynomial. 
Such an $e(z)$ can be realized as a vertex operator. 
An explicit formula is given by $E_{l+1}(z)$ in (\ref{defofE:k}).
For $k=2$, $e(z)$ satisfies 
$e(z)^3=0$ and $e(z)^2e(-z)=0$. 
Such an operator can be constructed as a sum $e(z)=a_1(z)+a_2(z)$, with
$[a_i(z),a_i(w)]=0$, $a_i(z)^2=0$ and $a_1(z)a_2(-z)=0$. 
Since the relations for $a_i(z)$ are quadratic, 
they can be realized as vertex operators. 

In general, for even $k=2s$, we set
\begin{eqnarray*}
&&e(z)=a_1(z)+\cdots+a_s(z)+b_1(z)+\cdots+b_s(z),
\\
&&[a_i(z),a_j(w)]=0,~~[b_i(z),b_j(w)]=0,~~[a_i(z),b_j(w)]=0,
\\
&&a_i(z)^2=0,~~b_i(z)^2=0, ~~a_i(z)b_i(-z)=0.
\end{eqnarray*}
For odd $k=2s+1$, 
\begin{eqnarray*}
&&e(z)=a_1(z)+\cdots+a_s(z)+c(z)+b_1(z)+\cdots+b_s(z),
\end{eqnarray*}
where $a_i(z),b_i(z)$ satisfy the same relations as above, and
\begin{eqnarray*}
&&[c(z),a_i(w)]=0,~~[c(z),b_i(w)]=0,~~c(z)^2=0,~~
c(z)c(-z)=0.
\end{eqnarray*}
Such $a_i(z),b_i(z),c(z)$ can be constructed as vertex operators,
and these operators are a part of a vertex operator realization of the 
entire algebra $\widehat{\mathfrak{sl}}_3$. 
The $\varepsilon$-deformation for even $k$ is given by
\begin{eqnarray*}
e_\varepsilon(z)=a_1(z)+\varepsilon a_2(z)+\cdots+\varepsilon^{s-1}a_s(z)
+\varepsilon^sb_1(z)+\cdots+\varepsilon^{2s-1}b_s(z),
\end{eqnarray*}
and similarly for odd $k$. 

The plan of this paper is as follows. 
Throughout this paper we consider $(k, 2)$ or $(k, 3)$-admissible configurations 
with the initial condition $a_{0} \le b_{0}$, see (\ref{boundary-admissible}). 
This corresponds to the case of $\Lambda=(k-b_{0})\Lambda_{0}+b_{0}\Lambda_{1}$ 
for $r=2$, and $\Lambda=b_{0}\Lambda_{1}+(k-b_{0})\Lambda_{2}$ for $r=3$. 
Here $\Lambda_{i}$'s are the fundamental weights of $\widehat{\mathfrak{sl}}_{r}$. 
As we mentioned above, the case $r=2$ has been studied in \cite{FS}. 
Nevertheless we give here the details in order to 
illustrate the method of vertex operators. 
The fermionic formulas for $r=3$ are new. 
First we introduce the algebra $E_{\Lambda}^{(k, r)}$ in Section 2. 
{}From Section 3 to Section 7 we consider the case of $r=2$. 
In Section 3 the dual space $(E_{\Lambda}^{(k, 2)})^{*}$ is realized 
as the space of functions $F^{(k, 2)}$. 
In order to calculate the character of $F^{(k, 2)}$ we define certain filtration 
$\{ \Gamma_{\lambda} \}$ on $F^{(k, 2)}$ in Section 4. 
Each component of the associated graded space determined by this filtration 
is embedded into a space of functions $G_{\lambda}^{(2)}\mathcal{S}_{\lambda}$,
see Proposition \ref{image:r=2}. 
We will prove that this embedding is surjective by using vertex operators. 
We summarize some properties of vertex operators constructed 
with $k$-dimensional bosons in Section 5. 
In Section 6 we give the current $e_{21}(z)$ using the vertex operator 
and prove that the dual space $\widetilde{W}_{0}^{*}$ is isomorphic to 
the space $G_{\lambda}^{(2)}\mathcal{S}_{\lambda}$. 
This implies surjectivity of the embedding and 
we get the fermionic formula for the character of $E_{\Lambda}^{(k, 2)}$, 
which is given in Section 7. 
In Section 8 we apply the argument above to the case of $r=3$. 
We use the two currents $e_{21}(z)$ and $e_{31}(z)$ and obtain the fermionic formula. 
As mentioned before we can construct a little different realization 
of the representation of $\widehat{\mathfrak{sl}}_{3}$ by using 
the mixed current $e(z)=e_{21}(z^{2})+z e_{31}(z^{2})$. 
This representation is of highest weight 
$\Lambda=[\frac{k+1}{2}]\Lambda_{1}+[\frac{k}{2}]\Lambda_{2}$. 
In Section 9 we obtain another type of fermionic formula in this special case 
using the current $e(z)$. 
This fermionic formula is the one obtained from a combinatorial point of view 
in \cite{combinatorial}. 
We give additional results and 
discuss some remaining problems in Section 10.

\setcounter{section}{1} 
\setcounter{equation}{0} 

\section{Preliminaries} 

\subsection{A polynomial algebra $E_{\Lambda}^{(k, r)}$} 

Let $r$ be a positive integer. 
Consider the polynomial ring 
\bea 
E^{(r)}:=\C[ e_{1}[-n], e_{2}[-n], \cdots , e_{r-1}[-n] ; n \ge 0]. 
\ena
We define formal power series $e_{j}(z)$ in $z$ by 
\bea 
e_{j}(z):=\sum_{n=0}^{\infty}e_{j}[-n]z^{n} , \quad (j=1, \ldots , r-1).  
\label{def-of-current} 
\ena 

Denote by $\{\Lambda_{i}\}_{i=0}^{r-1}$ the set of the fundamental weights of 
$\widehat{\mathfrak{sl}}_{r}$. 
Let $k$ be a positive integer and 
$\mathbf{b}=(b_{0}, \ldots , b_{r-2})$ a vector with non-negative integer entries such that 
\bea 
0 \le b_{0} \le \cdots \le b_{r-2} \le k. 
\ena 
We set the dominant integral weight $\Lambda$ of level $k$ by 
\bea 
\Lambda=(k-b_{r-2})\Lambda_{0}+b_{0}\Lambda_{1}+(b_{1}-b_{0})\Lambda_{2}+ \cdots 
         +(b_{r-2}-b_{r-3})\Lambda_{r-1}. 
\label{def-of-Lambda} 
\ena 

Denote by $J_{\Lambda}^{(k, r)}$ 
the ideal of $E^{(r)}$ generated by the elements 
\bea 
e_{1}[0]^{c_{0}} \cdots e_{i+1}[0]^{c_{i}}, \quad (i=0, \ldots , r-2), 
\label{boundary-vanish}
\ena 
where $c_{i}$'s are non-negative integers such that 
\bea 
c_{0}+\cdots +c_{i} > b_{i}, 
\label{boundary-cond} 
\ena 
and all the coefficients of the power series in the following form: 
\bea 
e_{1}(z)^{p_{1}} \cdots e_{r-1}(z)^{p_{r-1}},  
\ena 
where $p_{1}, \ldots , p_{r-1}$ are non-negative integers satisfying 
\bea 
p_{1}+ \cdots +p_{r-1}=k+1.  
\ena  

Set
\bea 
E_{\Lambda}^{(k, r)}:= E^{(r)}/J_{\Lambda}^{(k, r)}. 
\ena 
Now we give a basis of the vector space $E_{\Lambda}^{(k, r)}$. 
For $e \in E^{(r)}$ we denote by $\overline{e} \in E_{\Lambda}^{(k, r)}$ 
the image of $e$ by the projection 
$E^{(r)} \twoheadrightarrow E_{\Lambda}^{(k, r)}$. 
Let $\mathbf{a}=(a_{i})_{i=0}^{\infty}$ be a sequence of 
non-negative integers with finitely many non-zero entries. 
We define $e(\mathbf{a}) \in E^{(r)}$ by 
\bea 
e(\mathbf{a})&:=&\prod_{n \ge 0}\prod_{i=1}^{r-1} e_{i}[-n]^{a_{(r-1)n+i-1}} \nn \\ 
&=& \cdots e_{r-1}[-1]^{a_{2r-3}} \cdots e_{1}[-1]^{a_{r-1}} 
    e_{r-1}[0]^{a_{r-2}} \cdots e_{1}[0]^{a_{0}}.  
\ena 

A sequence $\mathbf{a}=(a_{i})_{i=0}^{\infty}$ of integers 
with finitely many non-zero entries is called $(k, r)$-{\it admissible} if 
\bea 
0 \le a_{i} \le k, \quad 
a_{i}+\cdots +a_{i+r-1} \le k 
\label{adm-cond} 
\ena 
for all $i \ge 0$. 
Denote by $C_{\mathbf{b}}^{(k, r)}$ the set of 
all $(k, r)$-admissible sequences such that 
\bea 
a_{0}\le b_{0}, \, a_{0}+a_{1} \le b_{1}, \,  \ldots , 
a_{0}+\cdots + a_{r-2} \le b_{r-2}. 
\label{boundary-admissible} 
\ena 

\begin{prop}\label{W-basis}
The set 
\bea 
\{ \overline{e(\mathbf{a})}; \mathbf{a} \in C^{(k, r)}_{\mathbf{b}}\} 
\label{spanning-set} 
\ena 
is a basis of $E_{\Lambda}^{(k, r)}$.
\end{prop} 

This proposition is a special case of the result by Primc \cite{P} 
which we will explain below. 
%
Now set $e_{i}[n]=e_{i+1, 1}[n] \in \widehat{\mathfrak{sl}}_{3}$. 
Then the elements (\ref{boundary-vanish}) satisfy 
\bea 
e_{1}[0]^{c_{0}} \cdots e_{i+1}[0]^{c_{i}}v_{\Lambda}=0
\label{e-boundary} 
\ena 
for non-negative integers $\{c_{i}\}$ satisfying (\ref{boundary-cond})
and the formal power series (\ref{def-of-current}) satisfy 
\bea 
e_{1}(z)^{p_{1}} \cdots e_{r-1}(z)^{p_{r-1}}=0
\label{integrability}
\ena 
on $L(\Lambda)$ for non-negative integers $p_{1}, \ldots , p_{r-1}$ 
such that $p_{1}+ \cdots +p_{r-1}=k+1.$ 
Hence the map 
\bea 
E_{\Lambda}^{(k, r)} \ni \overline{e} \,\,  \mapsto \,\,  
\overline{e}v_{\Lambda} \in W(\Lambda)
\label{iso-W} 
\ena 
is well-defined. Here $W(\Lambda)$ is the subspace defined by (\ref{principal-space}). 
This map is also surjective. 

In \cite{P}, Primc constructed a basis of $W(\Lambda)$. 
For $\mathbf{a}=(a_{i})_{i=0}^{\infty} \in C^{(k, r)}$, 
define the vector $M(\mathbf{a})$ of $W(\Lambda)$ by
\bea 
M(\mathbf{a}):=\overline{e(\mathbf{a})}v_{\Lambda}. 
\ena 
\begin{thm}\cite{P}  
Let $\Lambda$ be the dominant integral weight given by (\ref{def-of-Lambda}). 
Then the set 
\bea 
\mathcal{M}(\Lambda):=\{M(\mathbf{a}) ; \mathbf{a} \in C^{(k, r)}_{\mathbf{b}}\} 
\ena 
constitutes a basis of $W(\Lambda)$. 
\end{thm} 

{}From this theorem, the map (\ref{iso-W}) is injective and 
this implies Proposition \ref{W-basis}.

\subsection{Characters of $(k, r)$-admissible configurations} 

Now we introduce two kinds of degrees on $E^{(r)}$. 
First we define the {\it $q$-degree} by 
\bea 
{\rm deg}_{q}e_{i}[-n]:=(r-1)n+i-1. 
\ena 
Next define the {\it $z$-degree} by 
\bea 
{\rm deg}_{z}e_{i}[-n]:=1 
\ena 
for all $i=1, \ldots , r-1$ and $n \ge 0$. 

Note that the ideal $J_{\Lambda}^{(k, r)}$ is generated by 
homogeneous elements with respect to both of the degrees. 
Hence $E_{\Lambda}^{(k, r)}$ is a graded vector space 
with ${\rm deg}_{q}$ and ${\rm deg}_{z}$. 

Denote by $E_{\Lambda; i, j}^{(k, r)}$ the subspace 
spanned by homogeneous elements 
of $q$-degree $i$ and $z$-degree $j$. 
Consider the character 
\bea 
\chi_{E_{\Lambda}^{(k, r)}}(q, z):=
\sum_{i, j \ge 0}(\dim{E_{\Lambda; i, j}^{(k, r)}})q^{i}z^{j}. 
\ena 
{}From Proposition \ref{W-basis}, we have 
\bea 
\chi_{E_{\Lambda}^{(k, r)}}(q, z)=
\sum_{\mathbf{a}\in C_{\mathbf{b}}^{(k, r)}}
q^{\sum_{j \ge 0}ja_{j}}z^{\sum_{j \ge 0}a_{j}}. 
\ena 
This is nothing but the {\it character of $(k, r)$-configurations}
$\chi_{k, r; \mathbf{b}}(q, z)$ \cite{bosonic1}. 

In the following we give fermionic fomulas for the characters $\chi_{k, r; \mathbf{b}}$ 
in the two cases: 
\bea
{\rm (I)} && \!\!\!\!
r=2, (\Lambda=(k-b_{0})\Lambda_{0}+b_{0}\Lambda_{1}), \nn \\ 
{\rm (II)} && \!\!\!\!
r=3, b_{1}=k, \, (\Lambda=b_{0}\Lambda_{1}+(k-b_{0})\Lambda_{2}). \nn 
\ena 
In other words we consider $(k, 2)$ or $(k, 3)$-admissible configurations 
with the initial condition $a_{0} \le b_{0}$. 


\setcounter{section}{2} 
\setcounter{equation}{0} 

\section{Functional realization} 

{}From this section to Section \ref{end-of-r=2}, 
we consider $(k, 2)$-admissible configurations. 
In the following we fix $\Lambda=(k-b_{0})\Lambda_{0}+b_{0}\Lambda_{1}$ and 
abbreviate $E_{\Lambda}^{(k, 2)}$ and $J_{\Lambda}^{(k, 2)}$ 
to $E^{(k, 2)}$ and $J^{(k, 2)}$, respectively. 

Denote by $F_{n}$ the space of symmetric polynomials with $n$ variables: 
\bea 
F_{n}:=\C[x_{1}, \ldots , x_{n}]^{\mathfrak{S}_{n}}. 
\ena 
Let $E_{n}^{(2)}$ be the graded component of $E^{(2)}$ with $z$-degree $n$. 

We introduce a pairing 
\bea 
\langle \cdot , \cdot \rangle : 
E_{n}^{(2)} \otimes F_{n} \longrightarrow \C  
\ena 
as follows. 
Set $e(z):=e^{(1)}(z)$. 
Then we define the pairing by 
\bea 
\langle e(z_{1})\cdots e(z_{n}) , f(x_{1}, \ldots , x_{n}) \rangle 
:=f(z_{1}, \ldots , z_{n}). 
\label{def-pairing} 
\ena 

It is easy to see that the pairing $\langle \cdot , \cdot \rangle$ 
is a bilinear non-degenerate pairing. 
Moreover, it respects the grading on $E_{n}^{(2)}$ defined by the $q$-degree 
and the one on $F_{n}$ defined by the usual degree: $\deg x_{i}=1$. 

Denote by $J_{n}^{(k, 2)}$ the graded component of $J^{(k, 2)}$ with $z$-degree $n$. 

\begin{prop}\label{perp-space} 
The orthogonal complement $F_{n}^{(k, 2)}:=(J_{n}^{(k, 2)})^{\perp} \subset F_{n}$ 
is given as follows: 
\bea 
F_{n}^{(k, 2)}:=\left\{ f \in F_{n}; 
f(x_{1}, \ldots , x_{n})=0 \,\, {\rm if} 
\begin{array}{l} 
\,\, x_{1}=\cdots =x_{k+1} \,\, {\rm or}\\ 
\,\, x_{1}=\cdots =x_{b_{0}+1}=0. 
\end{array} \right\}. 
\ena 
\end{prop} 

\begin{proof} 
{}From the conditions (\ref{integrability}) and (\ref{e-boundary}) 
we have $e(z)^{k+1}=0$ and $e(0)^{b_{0}+1}=0$. 
Note that 
\bea 
\langle \textstyle{e(z)^{k+1} \prod_{j=k+2}^{n}e(z_{j})}, f(x_{1}, \ldots , x_{n}) \rangle 
=f(z, \ldots , z, z_{k+2}, \ldots , z_{n})
\label{pf-of-3.1-1} 
\ena 
and 
\bea 
\langle \textstyle{e(0)^{b_{0}+1} \prod_{j=b_{0}+2}^{n}e(z_{j})}, f(x_{1}, \ldots , x_{n}) \rangle 
=f(0, \ldots , 0, z_{b_{0}+2}, \ldots , z_{n}).
\label{pf-of-3.1-2} 
\ena 
Both of (\ref{pf-of-3.1-1}) and (\ref{pf-of-3.1-2}) equal zero 
if and only if $f \in F_{n}^{(k, 2)}$. 
\end{proof} 

Note that the graded components $E^{(2)}_{n}$ and $F_{n}$ are finite-dimensional 
and the pairing respects the grading. 
Therefore $(F_{n}^{(k, 2)})^{\perp}=(J_{n}^{(k, 2)})^{\perp\perp}=J_{n}^{(k, 2)}$  
and we obtained the following. 
\begin{prop} 
The pairing $\langle \cdot , \cdot \rangle$ induces 
a well-defined non-degenerate bilinear pairing of graded spaces 
\bea 
\langle \cdot , \cdot \rangle : 
E_{n}^{(k, 2)} \otimes F_{n}^{(k, 2)} \longrightarrow \C, 
\ena 
where $E_{n}^{(k, 2)}$ is the graded component of $E^{(k, 2)}$ with $z$-degree $n$. 
\end{prop} 

Hence the character $\chi_{k, 2; b_{0}}(q, z)$ 
is represented in terms of the character of $F_{n}^{(k, 2)}$ as follows. 
The character $\mathop{\rm ch}F_{n}^{(k, 2)}(q)$ is defined by 
\bea 
\mathop{\rm ch}F_{n}^{(k, 2)}(q):=\sum_{m=0}^{\infty} q^{m} \dim{( F_{n}^{(k, 2)})_{m} }, 
\ena 
where $( F_{n}^{(k, 2)})_{m}$ is the graded component of degree $m$. 
Then we get 
\begin{cor}\label{character-functional}
\bea 
\chi_{k, 2; b_{0}}(q, z)=\sum_{n=0}^{\infty}z^{n} \mathop{\rm ch}F_{n}^{(k, 2)}(q). 
\ena 
\end{cor}

\setcounter{section}{3} 
\setcounter{equation}{0} 

\section{Gordon filtration}\label{sec:2}

Let $k \in \Z_{\ge 0}$ and $n \in \Z_{\ge k}$. 
Let $\lambda$ be a level-$k$ restricted partition of $n$, that is 
\bea 
\lambda=(1^{m_{1}}, 2^{m_{2}}, \ldots, k^{m_{k}}), \quad 
\sum_{a=1}^{k}am_{a}=n. 
\ena 
Denote by $m_{a}(\lambda)$ the number of rows of length $a$ 
in the partition (or Young diagram) $\lambda$. 
Set $\mathbf{m}(\lambda):=(m_{1}(\lambda), \ldots , m_{k}(\lambda))$. 

For a sequence of non-negative integers 
$\mathbf{m}=(m_{1}, \ldots , m_{r})$, 
we define the space of functions $\mathcal{S}_{\mathbf{m}}$ by 
\bea 
\mathcal{S}_{\mathbf{m}}:= 
\C[x_{1}^{(1)}, \ldots , x_{m_{1}}^{(1)}]^{\mathfrak{S}_{m_{1}}} \otimes 
\cdots \otimes 
\C[x_{1}^{(r)}, \ldots , x_{m_{r}}^{(r)}]^{\mathfrak{S}_{m_{r}}}. 
\label{defofS_{m}} 
\ena 
In particular, for a level-$k$ restricted partition $\lambda$ of $n$, 
we abbreviate $\mathcal{S}_{\mathbf{m}(\lambda)}$ to $\mathcal{S}_{\lambda}$. 
Now we define a map 
\bea 
\varphi_{\lambda}: 
\C[x_{1}, \ldots , x_{n}]^{\mathfrak{S}_{n}} \longrightarrow
\mathcal{S}_{\lambda} 
\label{varphi:r=2} 
\ena 
as follows. 
Fix a numbering from $1$ to $n$
of the set of indices 
\bea 
\{ (a, i, j) ; 1 \le a \le k, 1 \le i \le m_{a}(\lambda), 1 \le j \le a\}. 
\ena 
We set $\varphi(x_{m}) := x_{i}^{(a)}$ where 
$(a, i, j)$ is the $m$-th index in this numbering. 
Then the map $\varphi_{\lambda}$ is defined by 
\bea 
\varphi_{\lambda}(f(x_{1}, \ldots , x_{n})):= 
f(\varphi(x_{1}), \ldots , \varphi(x_{n})). 
\ena 
Since $f$ is symmetric, this map does not depend on the numbering. 

Introduce the lexicographical order on partitions of $n$ by 
\bea 
\lambda \succ \mu \Longleftrightarrow 
\lambda_{j}=\mu_{j} \, ( j<p) \,\, {\rm and} \,\, 
 \lambda_{p}>\mu_{p} \,\, ,{\rm for \,\, some} \, p. 
\ena 
We define the subspaces of $F_{n}^{(k, 2)}$ by 
\bea 
\mathcal{F}_{\lambda} &:=& \mathop{\rm Ker} \varphi_{\lambda} \cap F_{n}^{(k, 2)}, \\ 
\Gamma_{\lambda} &:=& \cap_{\nu \succ \lambda} \mathcal{F}_{\nu}, \\
\Gamma'_{\lambda} &:=& \Gamma_{\lambda} \cap \mathop{\rm Ker} \varphi_{\lambda}. 
\ena 
The subspaces $\Gamma_{\lambda}$ give a filtration of $F_{n}^{(k, 2)}$ 
and we have 
\bea 
\mathop{\rm ch}F_{n}^{(k, 2)}=\sum_{\lambda}\mathop{\rm ch}
\left( \Gamma_{\lambda}/\Gamma'_{\lambda} \right), 
\label{fermionic-decomp} 
\ena 
where the right hand side is the summation 
over all level-$k$ restricted partitions of $n$. 

For an integer $s$ we set $(s)_{+}:={\rm max}(s, 0)$. 

\begin{prop}\label{image}
Let $\lambda$ be a level-$k$ restricted partition of $n$. 
The image of the map $\varphi_{\lambda}|_{\Gamma_{\lambda}}$ is contained 
in the principal ideal $G_{\lambda}^{(2)}\mathcal{S}_{\lambda}$, 
where the function $G_{\lambda}^{(2)}$ is defined by 
\bea 
G_{\lambda}^{(2)}:= 
\prod_{a=1}^{k}\prod_{j}(x_{j}^{(a)})^{(a-b_{0})_{+}} \!\!\!
\prod_{1 \le a<b \le k}\prod_{i, j}
(x_{i}^{(a)}-x_{j}^{(b)})^{2a} 
\prod_{a=1}^{k}\prod_{i<j}
(x_{i}^{(a)}-x_{j}^{(a)})^{2a}. 
\label{def-of-G:r=2}
\ena 
Hence the map $\varphi_{\lambda}|_{\Gamma_{\lambda}}$ induces the embedding 
of the subquotient $\Gamma_{\lambda}/\Gamma'_{\lambda}$ 
into the principal ideal $G_{\lambda}^{(2)}\mathcal{S}_{\lambda}$. 
\end{prop}

\begin{proof} 
Similar to the proof of Lemma 3.4.1 and Lemma 3.4.3 in \cite{coinv3}. 
\end{proof} 

Our goal is to prove that the image of $\varphi_{\lambda}|_{\Gamma_{\lambda}}$
is equal to $G^{(2)}_{\lambda}\mathcal{S}_{\lambda}$.

\setcounter{section}{4}
\setcounter{equation}{0}

\section{Vertex operators}\label{sec:3}

\subsection{Definitions}
Let $N$ be a positive integer. 
We fix a non-degenerate symmetric bilinear form $\langle \cdot , \cdot \rangle$ 
on the $N$-dimensional $\C$-vector space $\C^{N}$. 

We denote by $\widehat{\mathcal{H}}_{N}$ the Heisenberg algebra with unit $1$
generated by the elements $a_{m}(\alpha)$ and $e^{Q(\alpha)}$ 
$(m \in \Z, \alpha \in \C^{N})$ 
satisfying the relations
\bea
&& 
[ a_{m}(\alpha), a_{n}(\beta) ]=
m\langle \alpha, \beta \rangle\delta_{m+n, 0}, \\ 
&& 
[ a_{m}(\alpha), e^{Q(\beta)} ]= 
\delta_{m, 0}\langle \alpha , \beta \rangle e^{Q(\beta)}, \quad 
e^{Q(\alpha)}e^{Q(\beta)}=e^{Q(\alpha+\beta)}. 
\ena 
Here the generators $a_{m}(\alpha)$ are linear on $\alpha$. 

We define the Fock space $\mathcal{F}$ by 
\bea 
\mathcal{F}:=\C[ a_{-m}(\alpha) ; m>0, \alpha \in \C^{N}] \otimes 
             \C[ e^{Q(\beta)} ; \beta \in \C^{N}]. 
\ena 
The algebra $\widehat{\mathcal{H}}_{N}$ acts on $\mathcal{F}$ 
as follows: 
\bea 
&& 
a_{m}(\alpha) (f \otimes e^{Q(\beta)})= \left\{ 
\begin{array}{lc} 
(a_{m}(\alpha)f) \otimes e^{Q(\beta)}, & (m<0), \\ 
{}[a_{m}(\alpha), f] \otimes e^{Q(\beta)}, & (m > 0), \\ 
\langle \alpha, \beta \rangle f \otimes e^{Q(\beta)}, & (m=0), 
\end{array} \right. \\ 
&& 
e^{Q(\alpha)} (f \otimes e^{Q(\beta)})= 
f \otimes e^{Q(\alpha+\beta)}, 
\ena 
where $f \in \C[ a_{-m}(\alpha) ; m>0, \alpha \in \C^{N}]$.

Let $\boldsymbol\alpha=(\{\alpha_{m}\}_{m \in \Z}, \alpha^{0})$ 
be a sequence of vectors in $\C^{N}$. 
The vertex operator $X_{\boldsymbol\alpha}(z)$ is defined by 
\bea 
X_{\boldsymbol\alpha}(z):=
\exp{\Bigl( \sum_{m>0}\frac{a_{-m}(\alpha_{-m})}{m}z^{m} \Bigr)} 
\exp{\Bigl( -\sum_{m>0}\frac{a_{m}(\alpha_{m})}{m}z^{-m} \Bigr)} 
e^{Q(\alpha^{0})}z^{a_{0}(\alpha_{0})}. 
\ena 
Introduce the normal ordering $: \cdot :$ on $\widehat{\mathcal{H}}_{N}$: 
\bea 
&& 
:a_{m}(\alpha)a_{n}(\beta)\!: \, = \left\{
\begin{array}{lc}  
a_{m}(\alpha)a_{n}(\beta), & (m<0), \\ 
a_{n}(\beta)a_{m}(\alpha), & (m>0), 
\end{array} \right. \\
&& 
:a_{0}(\alpha) e^{Q(\beta)}\!: \, = \, 
:e^{Q(\beta)} a_{0}(\alpha)\!: \, = 
e^{Q(\beta)} a_{0}(\alpha). 
\ena 
Then we have 
\bea 
X_{\boldsymbol\alpha}(z)X_{\boldsymbol\beta}(w)= 
g(z, w; \boldsymbol\alpha, \boldsymbol\beta) 
:X_{\boldsymbol\alpha}(z)X_{\boldsymbol\beta}(w)\!:, 
\ena 
where
\bea 
g(z, w ; \boldsymbol\alpha, \boldsymbol\beta ):= 
z^{\langle \alpha_{0}, \beta^{0} \rangle} 
\exp{\Bigl( -\sum_{m>0}\frac{\langle \alpha_{m}, \beta_{-m} \rangle}{m} 
        \left( \frac{w}{z} \right)^{m} \Bigr)} 
\ena 
for $\boldsymbol\alpha=(\{\alpha_{m}\}, \alpha^{0})$ and 
$\boldsymbol\beta=(\{\beta_{m}\}, \beta^{0})$.

\subsection{Matrix elements}
Set 
\bea 
\ket{\beta}:= 1 \otimes e^{Q(\beta)} \in \mathcal{F}. 
\ena 
Note that 
\bea 
a_{m}(\alpha)\ket{\beta}=0, \quad {\rm if} \quad m>0. 
\ena 
Let $\bra{\beta} \in \mathcal{F}^{*}$ be the dual vector defined by 
\bea 
\bra{\beta}(f \otimes e^{Q(\gamma)})= \left\{ 
\begin{array}{ll} 
c, & {\rm if}\,\, f=c \in \C \,\, {\rm and} \,\, \beta=\gamma, \\ 
0, & {\rm otherwise}. 
\end{array} \right. 
\ena 

Denote by $\widehat{\mathcal{H}}_{N}^{+}$ 
the commutative subalgebra of $\widehat{\mathcal{H}}_{N}$ 
generated by the generators $a_{m}(\alpha), (m>0, \alpha \in \C^{N})$ 
and $1$. 
Consider the matrix element 
\bea 
\bra{\beta'} h
X_{\boldsymbol\alpha_{1}}(x_{1}) \cdots 
X_{\boldsymbol\alpha_{n}}(x_{n}) \ket{\beta}, 
\qquad h \in \widehat{\mathcal{H}}_{N}^{+},
\ena 
for $\boldsymbol\alpha_{a}=(\{\alpha_{a, m}\}, \alpha_{a}^{0}), 
(a=1, \ldots , n)$. 
{}From the definiton of $X_{\boldsymbol\alpha}(z)$, 
it is easy to see that 
\bea 
\bra{\beta} h
X_{\boldsymbol\alpha_{1}}(x_{1}) \cdots 
X_{\boldsymbol\alpha_{n}}(x_{n}) \ket{0}=0 \quad 
{\rm unless} \quad \beta'-\beta=\sum_{a=1}^{n}\alpha_{a}^{0}. 
\ena 

\begin{thm}\label{matrixelt} 
Let $\boldsymbol\alpha_{a}=(\{\alpha_{a, m}\}, \alpha_{a}^{0}), \, 
 (a=1, \ldots , N)$ be sequences of vectors in $\C^{N}$ 
and $\mathbf{m}=(m_{1}, \ldots , m_{N}), (\forall m_{j} \ge 0)$ 
a sequence of non-negative integers. 
Denote by 
$S_{\mathbf{m}}(\boldsymbol\alpha_{1}, \cdots , \boldsymbol\alpha_{N}; \beta)$ 
the set of functions given by 
\bea 
\{\bra{\beta+\alpha^{*}} h 
X_{\boldsymbol\alpha_{1}}(x^{(1)}_{1}) \cdots 
X_{\boldsymbol\alpha_{1}}(x^{(1)}_{m_{1}}) \cdots 
X_{\boldsymbol\alpha_{N}}(x^{(N)}_{1}) \cdots 
X_{\boldsymbol\alpha_{N}}(x^{(N)}_{m_{N}}) \ket{\beta}; 
h \in \widehat{\mathcal{H}}_{N}^{+}\}, 
\label{defofV} 
\ena 
where $\alpha^{*}:=\textstyle{\sum_{a}m_{a}\alpha_{a}^{0}}$ 
and $\beta \in \C^{N}$. 

Suppose that 
the vectors $\alpha_{a, -m} \, (a=1, \ldots , N)$ are 
linearly independent for each $m>0$. 
Then we have 
\bea 
S_{\mathbf{m}}(\boldsymbol\alpha_{1}, \cdots , \boldsymbol\alpha_{N}; \beta)
&=& 
\prod_{a=1}^{N}\prod_{j}(x_{j}^{(a)})^{\langle \alpha_{a}^{0}, \beta \rangle} 
\!\!\! 
\prod_{1 \le a<b \le N}\prod_{i, j}
g(x_{i}^{(a)}, x_{j}^{(b)} ; \boldsymbol\alpha_{a}, \boldsymbol\alpha_{b}) \nn \\ 
&\times& 
\prod_{a=1}^{N}\prod_{i<j}
g(x_{i}^{(a)}, x_{j}^{(a)} ; \boldsymbol\alpha_{a}, \boldsymbol\alpha_{a}) 
\cdot \mathcal{S}_{\mathbf{m}}
\ena 
for any $\mathbf{m}$. 
\end{thm}

\begin{proof} 
Fix $\mathbf{m}=(m_{1}, \ldots , m_{N})$. 
For $h \in \widehat{\mathcal{H}}_{N}^{+}$, we set 
\bea 
F(h):=
\bra{\beta+\alpha^{*}} h 
X_{\boldsymbol\alpha_{1}}(x^{(1)}_{1}) \cdots 
X_{\boldsymbol\alpha_{1}}(x^{(1)}_{m_{1}}) \cdots 
X_{\boldsymbol\alpha_{N}}(x^{(N)}_{1}) \cdots 
X_{\boldsymbol\alpha_{N}}(x^{(N)}_{m_{N}}) \ket{\beta}. 
\ena 
Then it is easy to see that 
\bea 
F(1)&=& 
\prod_{a=1}^{N}\prod_{j}(x_{j}^{(a)})^{\langle \alpha_{a}^{0}, \beta \rangle} 
\!\!\! 
\prod_{1 \le a<b \le N}\prod_{i, j}
g(x_{i}^{(a)}, x_{j}^{(b)} ; \boldsymbol\alpha_{a}, \boldsymbol\alpha_{b}) \nn \\
&\times& 
\prod_{a=1}^{N}\prod_{i<j}
g(x_{i}^{(a)}, x_{j}^{(a)} ; \boldsymbol\alpha_{a}, \boldsymbol\alpha_{a}). 
\ena 

For $r>0$ the vertex operator $X_{\boldsymbol\alpha}(z)$ satisfies 
\bea 
[a_{r}(\gamma), X_{\boldsymbol\alpha}(z)]=
\langle \gamma , \alpha_{-r} \rangle z^{r} 
X_{\boldsymbol\alpha}(z). 
\ena 
Hence we have 
\bea 
F(a_{r_{1}}(\gamma_{1}) \cdots a_{r_{l}}(\gamma_{l}))= 
\prod_{i=1}^{l} \left( 
\sum_{a=1}^{N}
\langle \gamma_{i}, \alpha_{a, -r_{i}} \rangle 
p_{r_{i}}^{(a)} \right) F(1), \quad 
(\forall r_{i}>0), 
\ena 
where $p_{r}^{(a)}$ is the $r$-th power sum of $x_{j}^{(a)}$'s, 
that is $p_{r}^{(a)}:=\sum_{j}(x_{j}^{(a)})^{r}$. 
Therefore, 
if the vectors $\alpha_{a, -r}, (a=1, \ldots , N)$ are linearly independent 
for each $r>0$, 
we can obtain any polynomial in $\mathcal{S}_{\mathbf{m}}$ 
as $F(h)/F(1)$ by taking a suitable $h \in \widehat{\mathcal{H}}_{N}^{+}$.  
\end{proof}

\setcounter{section}{5}
\setcounter{equation}{0}

\section{Construction of vertex operators}\label{sec:4}

Fix a basis $\{\epsilon_{a}\}_{a=1}^{k}$ of $\C^{k}$ satisfying 
\bea 
\langle \epsilon_{a}, \epsilon_{b} \rangle = 2\delta_{a, b}. 
\ena 
For $1 \le a \le k$, define a sequence of vectors 
$\boldsymbol\alpha_{a}=(\{\alpha_{a,m}\}, \alpha_{a}^{0})$ by 
\bea 
\alpha_{a, m}=\alpha^{0}_{a}=\epsilon_{a}, \quad (\forall m \in \Z). 
\ena 
Now we set 
\bea 
E_{a}(z):=X_{\boldsymbol\alpha_{a}}(z). 
\ena 
Then we have 
\bea 
E_{a}(z)E_{b}(w)= \left\{ 
\begin{array}{lc} 
:E_{a}(z)E_{b}(w)\!:, & a \not=b, \\ 
(z-w)^2 : E_{a}(z)E_{a}(w)\!:, & a=b. 
\end{array} \right.  
\ena 
In particular the operators $E_{a}(z)$ are commutative and satisfy 
$E_{a}(z)^{2}=0$. 

Set 
\bea 
E_{\varepsilon}(z):=\varepsilon_{1}E_{1}(z)+ \cdots +\varepsilon_{k}E_{k}(z). 
\ena 
Let $\lambda$ be a level-$k$ restricted partition of $n$ 
and $\lambda'=(\lambda_{1}', \ldots , \lambda_{k}')$ its conjugate (or transpose). 
Define the operator $E_{\lambda}$ by 
\bea 
E_{\lambda}(x_{1}, \ldots , x_{n}):=
\prod_{a=1}^{k}\frac{1}{\lambda_{a}'!} \Bigl.
\left(\frac{\partial}{\partial \varepsilon_{1}}\right)^{\lambda_{1}'} \cdots 
\left(\frac{\partial}{\partial \varepsilon_{k}}\right)^{\lambda_{k}'}
E_{\varepsilon}(x_{1}) \cdots E_{\varepsilon}(x_{n})
\Bigr|_{\forall \varepsilon_{a}=0}. 
\ena
In other words the operator $E_{\lambda}(x_{1}, \ldots , x_{n})$ 
is the symmetrization of 
\bea 
\prod_{a=1}^{k}(E_{a}(x_{n_{a-1}'+1}) \cdots E_{a}(x_{n_{a}'})), 
\label{symm-source} 
\ena  
where $n_{0}'=0$ and $n_{a}':=\sum_{j=1}^{a}\lambda_{j}'$. 

Set $\epsilon_{\lambda}^{*}:=\sum_{a}\lambda_{a}'\epsilon_{a}$. 
Note that 
\bea 
\bra{\beta'}h E_{\lambda}(x_{1}, \cdots , x_{n})\ket{\beta}=0 \,\, 
(\forall h \in {\widehat{\mathcal{H}}_{k}^{+}}), 
\quad 
{\rm unless} \,\,  \beta'-\beta=\epsilon_{\lambda}^{*}. 
\ena 
Consider the space of symmetric polynomials 
\bea 
U_{\lambda}:= \{
\bra{\beta_{0}+\epsilon_{\lambda}^{*}}h 
E_{\lambda}(x_{1}, \cdots , x_{n})\ket{\beta_{0}}; 
h \in \widehat{\mathcal{H}}_{k}^{+} \}, 
\ena 
where $\beta_{0}$ is the vector in $\C^{k}$ uniquely determined by 
\bea 
\langle \epsilon_{a}, \beta_{0} \rangle=\left\{ 
\begin{array}{ll} 
0, & {\rm if} \,\, 1 \le a \le b_{0} , \\ 
1, & {\rm if} \,\, a > b_{0}. 
\end{array} \right. 
\ena 

\begin{prop}\label{subset:r=2}
\bea 
U_{\lambda} \subset \Gamma_{\lambda}.
\ena 
\end{prop} 

\begin{proof} 
Set 
\bea 
F_{\lambda}(h; x_{1}, \ldots , x_{n}):=\bra{\beta_{0}+\epsilon_{\lambda}^{*}}h 
E_{\lambda}(x_{1}, \cdots , x_{n})\ket{\beta_{0}}. 
\ena 
It suffices to prove 
\bea 
\varphi_{\nu}(F_{\lambda}(h; x_{1}, \ldots , x_{n}))=0 
\label{pf-of-6.1-1} 
\ena 
for any $\nu \succ \lambda$ and 
\bea 
F_{\lambda}(h; 0, \ldots , 0, x_{b_{0}+2}, \ldots , x_{n})=0. 
\label{pf-of-6.1-2} 
\ena 

First we prove (\ref{pf-of-6.1-1}). 
Note that $E_{\lambda}$ is also the symmetrization of 
\bea 
\prod_{j=1}^{\lambda_{1}'}
\left( E_{1}(x_{n_{j-1}+1}) E_{2}(x_{n_{j-1}+2}) \cdots E_{\lambda_{j}}(x_{n_{j}}) \right), 
\label{factored-form} 
\ena 
where $n_{0}=0$ and $n_{j}:=\sum_{i=1}^{j}\lambda_{i}$. 
{}From this expression and the relation $E_{a}(z)^{2}=0$, 
it is easy to see $ \varphi_{\nu}(F_{\lambda})=0 $ for $\nu \succ \lambda$. 

Next we prove (\ref{pf-of-6.1-2}). 
Note that $E_{\lambda}$ is the symmetrization of (\ref{symm-source}).  
Consider the function 
\bea 
\bra{\beta_{0}+\epsilon_{\lambda}^{*}} h 
\prod_{a=1}^{k}(E_{a}(x_{n_{a-1}'+1}) \cdots E_{a}(x_{n_{a}'})) \ket{\beta_{0}}. 
\label{pf-of-6.1-3} 
\ena 
{}From Theorem \ref{matrixelt} this function (\ref{pf-of-6.1-3}) is a polynomial 
with the factor 
\bea 
\prod_{j=n_{b_{0}}'+1}^{n}x_{j} 
\prod_{1 \le a \le k \atop n_{a-1}' <i<j\le n_{a}'}(x_{i}-x_{j})^{2a}. 
\ena 
Hence (\ref{pf-of-6.1-3}) becomes zero 
if $b_{0}+1$ variables of $x_{1}, \ldots , x_{n}$ are equal to zero. 
This implies (\ref{pf-of-6.1-2}) 
because the function $F_{\lambda}$ is the symmetrization of (\ref{pf-of-6.1-3}).
\end{proof} 

\begin{prop}\label{image:r=2}
\bea 
\varphi_{\lambda}(U_{\lambda})=
S_{\mathbf{m}(\lambda)}(\boldsymbol\gamma_{1}, \cdots , \boldsymbol\gamma_{k}; \beta_{0}), 
\ena 
where $S_{\mathbf{m}(\lambda)}$ is the space defined in Theorem \ref{matrixelt} 
and the sequences of vectors 
$\boldsymbol\gamma_{a}=(\{\gamma_{a, m}\}, \gamma_{a}^{0}), (a=1, \ldots , k)$ 
are defined by 
\bea 
\textstyle{\gamma_{a, m}=\gamma_{a}^{0}=\sum_{j=1}^{a}\epsilon_{j}, \quad 
(\forall m \in \Z)}. 
\label{defofgamma:r=2}
\ena 
\end{prop} 

\begin{proof}
{}From the relation $E_{a}(z)^2=0$, we have 
\bea 
\varphi_{\lambda}(E_{\lambda}(x_{1}, \ldots , x_{n}))= 
z_{\lambda} 
\prod_{a=1}^{k} \prod_{j=1}^{m_{a}(\lambda)}
\left(  
E_{1}(x_{j}^{(a)}) \cdots E_{a}(x_{j}^{(a)}) \right),  
\ena 
where $z_{\lambda}$ is a constant defined by 
$z_{\lambda}:=\prod_{a=1}^{k}(a!)^{m_{a}(\lambda)}$. 
Moreover, we see that 
\bea 
E_{1}(x) \cdots E_{a}(x)= \,
:\!E_{1}(x) \cdots E_{a}(x)\!: \, =X_{\boldsymbol\gamma_{a}}(x). 
\ena 
This completes the proof. 
\end{proof} 

Note that the vectors $\gamma_{a, -m}, (a=1, \ldots , k)$ 
defined in (\ref{defofgamma:r=2}) are linearly independent for each $m>0$. 
It is easy to check that 
\bea 
g(z, w; \boldsymbol\gamma_{a}, \boldsymbol\gamma_{b})= 
(z-w)^{2{\rm min}(a, b)}, \quad 
\langle \gamma_{a}^{0}, \beta_{0} \rangle =(a-b_{0})_{+}. 
\ena 
Therefore, from Proposition \ref{image}, 
Theorem \ref{matrixelt} and Proposition \ref{image:r=2}, 
we see 
\begin{cor}\label{cor:r=2}
\bea 
\varphi_{\lambda}(\Gamma_{\lambda})= 
G_{\lambda}^{(2)} \mathcal{S}_{\lambda}. 
\ena  
\end{cor}

\noindent{\it Example.} 
Consider the case of $k=3, b_{0}=1$ and $n=3$. 
Then the Gordon filtration is 
\bea 
F_{3}^{(3, 2)}=\Gamma_{(2, 1)} \supset \Gamma_{(1, 1, 1)} \supset \{0\}, 
\ena 
where 
\bea 
F_{3}^{(3, 2)}&=&
\left\{ f(x_{1}, x_{2}, x_{3}) \in \C[x_{1}, x_{2}, x_{3}]^{\mathfrak{S}_{3}}; 
   f=0 \,\, {\rm if} 
\begin{array}{l} 
x_{1}=x_{2}=x_{3} \,\, {\rm or} \\ 
x_{1}=x_{2}=0 
\end{array} \right\}. \\
\Gamma_{(1, 1, 1)}&=&
\{ f(x_{1}, x_{2}, x_{3}) \in F_{3}^{(3, 2)}; f=0 \,\, {\rm if} \,\, x_{1}=x_{2} \} \\ 
&=& 
\{ f(x_{1}, x_{2}, x_{3}) \in \C[x_{1}, x_{2}, x_{3}]^{\mathfrak{S}_{3}}; 
   f=0 \,\, {\rm if} \,\, x_{1}=x_{2} \}.
\ena 
The map $\varphi_{\lambda}$ is defined by 
\bea 
&& 
\varphi_{(2, 1)}: f(x_{1}, x_{2}, x_{3}) \mapsto 
f(x^{(2)}_{1}, x^{(2)}_{1}, x^{(1)}_{1}) \in \C[x^{(1)}_{1}] \otimes \C[x^{(2)}_{1}], \\ 
&& 
\varphi_{(1, 1, 1)}: f(x_{1}, x_{2}, x_{3}) \mapsto 
f(x^{(1)}_{1}, x^{(1)}_{2}, x^{(1)}_{3}) \in 
\C[x^{(1)}_{1}, x^{(1)}_{2}, x^{(1)}_{3}]^{\mathfrak{S}_{3}}. 
\ena 
Corollary \ref{cor:r=2} shows that 
\bea 
&& 
\varphi_{(2, 1)}(\Gamma_{(2, 1)})=x_{1}^{(2)}(x^{(1)}_{1}-x^{(2)}_{1})^{2} 
\C[x^{(1)}_{1}] \otimes \C[x^{(2)}_{1}], \\ 
&& 
\varphi_{(1, 1, 1)}(\Gamma_{(1, 1, 1)})=
\prod_{1 \le i < j \le 3}\!\!\!(x^{(1)}_{i}-x^{(1)}_{j})\,
\C[x^{(1)}_{1}, x^{(1)}_{2}, x^{(1)}_{3}]^{\mathfrak{S}_{3}}. 
\ena

\setcounter{section}{6}
\setcounter{equation}{0}

\section{Fermionic formula} \label{end-of-r=2}

Recall Corollary \ref{character-functional}: we have 
\bea 
\chi_{k, 2; b_{0}}(q, z)=\sum_{n=0}^{\infty}z^{n} \mathop{\rm ch}F_{n}^{(k, 2)}(q). 
\ena 

Now let us write down the character of $F_{n}^{(k, 2)}$. 
{}From Proposition \ref{image} and Corollary \ref{cor:r=2}, 
we find 
\bea 
\mathop{\rm ch}(\Gamma_{\lambda}/\Gamma'_{\lambda})= 
\mathop{\rm ch}(G_{\lambda}^{(2)}\mathcal{S}_{\lambda}).
\ena 

It is easy to obtain the formula for ${\rm ch}(G_{\lambda}^{(2)}\mathcal{S}_{\lambda})$. 
Introduce the $k \times k$ matrix $A^{(2)}$ defined by 
\bea 
A^{(2)}=(A^{(2)}_{ab})_{1 \le a, b \le k}, \quad 
A^{(2)}_{ab}:=2{\rm min}(a, b). 
\label{gordon-mat:r=2} 
\ena 
Denote by $\mathbf{c}^{(2)}_{b_{0}}$ the vector defined by 
\bea 
\mathbf{c}^{(2)}_{b_{0}}:=(0, \ldots , 0, 1, 2, \ldots , k-b_{0}). 
\label{boundary-term:r=2} 
\ena 
Then we have 
\bea 
\mathop{\rm ch}(\Gamma_{\lambda}/\Gamma'_{\lambda})= 
\mathop{\rm ch}(G_{\lambda}^{(2)}\mathcal{S}_{\lambda})= 
\frac{q^{\frac{1}{2}({}^{t}\mathbf{m} A^{(2)} \mathbf{m}
          -({\rm diag}A^{(2)}) \cdot \mathbf{m}) 
          +{\mathbf{c}}^{(2)}_{b_{0}}\cdot\mathbf{m} }} 
     {(q)_{m_{1}(\lambda)} \cdots (q)_{m_{k}(\lambda)}}, 
\label{formula-GS} 
\ena 
where $\mathbf{m}={}^{t}\mathbf{m}(\lambda)={}^{t}(m_{1}(\lambda), \ldots , m_{k}(\lambda))$ 
and $(q)_{n}:=\prod_{j=1}^{n}(1-q^{j})$. 
Here the numerator in (\ref{formula-GS}) represents the degree of $G_{\lambda}^{(2)}$ 
and the part 
\bea 
\frac{1}{(q)_{m_{1}(\lambda)} \cdots (q)_{m_{k}(\lambda)}}
\ena 
is the character of $\mathcal{S}_{\lambda}$. 

Substituting (\ref{formula-GS}) into (\ref{fermionic-decomp}), 
we get the fermionic formula for $(k, 2)$-admissible configurations:
\begin{thm} 
\bea 
\chi_{2, r; b_{0}}(q, z)=\sum_{n=0}^{\infty} 
\sum_{m_{1}+2m_{2}+\cdots +km_{k}=n \atop m_{1}, \ldots , m_{k} \ge 0} 
\frac{q^{\frac{1}{2}({}^{t}\mathbf{m} A^{(2)} \mathbf{m}-({\rm diag}A^{(2)}) \cdot \mathbf{m}) 
           +\mathbf{c}^{(2)}_{b_{0}}\cdot\mathbf{m}}} 
     {(q)_{m_{1}} \cdots (q)_{m_{k}}} z^{n}, 
\ena 
where $A^{(2)}$ is the $k \times k$ matrix 
defined by (\ref{gordon-mat:r=2}), 
$\mathbf{c}^{(2)}_{b_{0}}$ is the vector defined by (\ref{boundary-term:r=2}) 
and $\mathbf{m}={}^{t}(m_{1}, \ldots , m_{k})$. 
\end{thm}

\setcounter{section}{7}
\setcounter{equation}{0}

\section{Fermionic formula for $\chi_{k, 3}$}\label{section:2k}

In this section we consider the case where $r=3$ and $\mathbf{b}=(b_{0}, k)$. 
We fix $\Lambda=b_{0}\Lambda_{1}+(k-b_{0})\Lambda_{2}$ 
and abbreviate $E_{\Lambda}^{(k, 3)}$ and $J_{\Lambda}^{(k, 3)}$ 
to $E^{(k, 3)}$ and $J^{(k, 3)}$, respectively. 

\subsection{Functional realization of $W^{(k, 3)}$} 
Consider the space of polynomials 
\bea 
F_{l_{1}, l_{2}}:=\C[x_{1}^{2}, \ldots , x_{l_{1}}^{2}]^{\mathfrak{S}_{l_{1}}} \otimes 
           \C[y_{1}^{2}, \ldots , y_{l_{2}}^{2}]^{\mathfrak{S}_{l_{2}}}
           \cdot {\textstyle \prod_{j=1}^{l_{2}}y_{j}}. 
\ena 
Let us introduce a pairing 
\bea 
\langle \cdot , \cdot \rangle : 
E^{(3)}_{n} \otimes (\bigoplus_{l_{1}+l_{2}=n \atop l_{1}, l_{2} \ge 0} F_{l_{1}, l_{2}}) 
\longrightarrow \C 
\ena 
as follows. 
Set 
\bea 
a(z):=e^{(1)}(z^{2}), \quad b(z):=ze^{(2)}(z^{2}). 
\ena 
Then we define the pairing by 
\bea 
&& 
\langle a(z_{1}) \cdots a(z_{l_{1}}) b(w_{1}) \cdots b(w_{l_{2}}), 
 f(x_{1}, \ldots , x_{m_{1}}; y_{1}, \ldots , y_{m_{2}}) \rangle \nn \\ 
&& \qquad {}:=
\delta_{l_{1}m_{1}}\delta_{l_{2}m_{2}}f(z_{1}, \ldots , z_{l_{1}}; w_{1}, \ldots , w_{l_{2}}) 
\label{another-pairing}
\ena 
for $f \in F_{m_{1}, m_{2}}$. 
This pairing is non-degenerate and 
respects the gradings on $E_{n}^{(3)}$ and $\oplus_{l_{1}+l_{2}=n}F_{l_{1}, l_{2}}$. 
Here the grading on $\oplus_{l_{1}+l_{2}=n}F_{l_{1}, l_{2}}$ is the usual one 
defined by $\deg x_{i}=1=\deg y_{i}$. 

Let us determine the orthogonal complement 
$F_{l_{1}, l_{2}}^{(k, 3)}:=(J_{n}^{(k, 3)})^{\perp} \cap F_{l_{1}, l_{2}}$ 
with respect to the pairing defined above. 
Denote by $I_{l_{1}, l_{2}}^{(k, 3)}$ the space of functions 
\bea 
g(x_{1}, \ldots , x_{l_{1}}; y_{1}, \ldots , y_{l_{2}}) 
\in 
\C[x_{1}, \ldots , x_{l_{1}}]^{\mathfrak{S}_{l_{1}}} \otimes 
\C[y_{1}, \ldots , y_{l_{2}}]^{\mathfrak{S}_{l_{2}}} 
\ena 
such that
\bea  
g=0 \quad {\rm if} && 
x_{1}= \cdots =x_{a}=y_{1}= \cdots =y_{b}, \,\,  
(a \ge 0, b \ge 0, a+b=k+1), \nn \\ 
&& {\rm or} \,\,\, 
x_{1}= \cdots =x_{b_{0}+1}=0. 
\label{realization:2k} 
\ena 
\begin{prop}\label{perp-space:another}
\bea 
F_{l_{1}, l_{2}}^{(k, 3)}=
\{ g(x_{1}^{2}, \ldots , x_{l_{1}}^{2}; y_{1}^{2}, \ldots , y_{l_{2}}^{2}) 
   \textstyle{\prod_{j=1}^{l_{2}}y_{j}}; g \in I_{l_{1}, l_{2}}^{(k, 3)} \}. 
\label{defofF:2k} 
\ena 
\end{prop} 
The proof is quite similar to that of Proposition \ref{perp-space}. 

{}From this proposition, we have 
\begin{prop} 
The pairing (\ref{another-pairing}) induces 
a well-defined non-degenerate bilinear pairing of the graded spaces 
\bea 
\langle \cdot , \cdot \rangle : 
E_{n}^{(k, 3)} \otimes 
(\bigoplus_{l_{1}+l_{2}=n \atop l_{1}, l_{2} \ge 0} F_{l_{1}, l_{2}}^{(k, 3)}) 
\longrightarrow \C. 
\ena  
\end{prop} 

Introduce the usual grading on $I_{l_{1}, l_{2}}^{(k, 3)}$ and 
denote by $\mathop{\rm ch}I_{l_{1}, l_{2}}^{(k, 3)}(q)$ the character of 
the graded space $I_{l_{1}, l_{2}}^{(k, 3)}$ with the formal variable $q$. 
{}From Proposition \ref{perp-space:another} we have 
\bea 
\mathop{\rm ch}F_{l_{1}, l_{2}}^{(k, 3)}(q)=
q^{l_{2}}\mathop{\rm ch}I_{l_{1}, l_{2}}^{(k, 3)}(q^{2}). 
\ena 
Hence the character of $(k, 3)$-admissible configurations is given as follows. 
\begin{cor}\label{fermionic:another}
\bea 
\chi_{k, 3; (b_{0}, k)}(q, z)=\sum_{n=0}^{\infty}\sum_{l_{1}+l_{2}=n \atop l_{1}, l_{2} \ge 0} 
z^{n} q^{l_{2}} \mathop{\rm ch}I_{l_{1}, l_{2}}^{(k, 3)}(q^{2}). 
\ena 
\end{cor}

\subsection{Gordon filtration} 

Let us introduce a filtration on $I_{l_{1}, l_{2}}^{(k, 3)}$. 
 
For a partition $\lambda$ of $n$, 
let us write clearly the variables in $\mathcal{S}_{\lambda}$ 
by $\mathcal{S}_{\lambda}=\mathcal{S}_{\lambda}(x)$. 
Let $\lambda$ and $\mu$ be level-$k$ restricted partitions of 
$l_{1}$ and $l_{2}$, respectively. 
We denote by $\varphi_{\lambda, \mu}$ 
the tensor product of $\varphi_{\lambda}$ and $\varphi_{\mu}$:
\bea 
\varphi_{\lambda, \mu}:=\varphi_{\lambda} \otimes \varphi_{\mu}: 
\C[x_{1}, \ldots , x_{l_{1}}]^{\mathfrak{S}_{l_{1}}} \otimes 
\C[y_{1}, \ldots , y_{l_{2}}]^{\mathfrak{S}_{l_{2}}} 
\longrightarrow 
\mathcal{S}_{\lambda}(x) \otimes \mathcal{S}_{\mu}(y). 
\ena 

We define the lexicographical order on pairs of partitions by 
\bea 
(\lambda^{(1)}, \mu^{(1)}) \succ (\lambda^{(2)}, \mu^{(2)}) \, 
\Longleftrightarrow \lambda^{(1)} \succ \lambda^{(2)},  \,\,
{\rm or} \,\, \lambda^{(1)}=\lambda^{(2)} \, {\rm and} \, \mu^{(1)} \succ \mu^{(2)}. 
\ena 
Now let us define the subspaces of $I_{l_{1}, l_{2}}^{(k, 3)}$ by 
\bea 
\mathcal{F}_{\lambda, \mu} &:=& 
\mathop{\rm Ker} \varphi_{\lambda, \mu} \cap I_{l_{1}, l_{2}}^{(k, 3)}, \\ 
\Gamma_{\lambda, \mu} &:=& 
\cap_{(\nu, \kappa) \succ (\lambda, \mu)} \mathcal{F}_{\nu, \kappa}, 
\label{gordon-filt:2k} \\
\Gamma'_{\lambda, \mu} &:=& 
\Gamma_{\lambda, \mu} \cap \mathop{\rm Ker} \varphi_{\lambda, \mu}. 
\ena 
The subspaces $\Gamma_{\lambda, \mu}$ give a filtration of $I_{l_{1}, l_{2}}^{(k, 3)}$ 
and we have 
\bea 
\mathop{\rm ch}I_{l_{1}, l_{2}}^{(k, 3)}=\sum_{(\lambda, \mu)}\mathop{\rm ch}
\left( \Gamma_{\lambda, \mu}/\Gamma'_{\lambda, \mu} \right). 
\ena 

In the same way as the proof of Proposition \ref{image}, 
we can show the following: 

\begin{prop}\label{image:r=3'}
Let $\lambda$ and $\mu$ be level-$k$ restricted partitions 
of $l_{1}$ and $l_{2}$, respectively. 
The image of the map $\varphi_{\lambda, \mu}|_{\Gamma_{\lambda, \mu}}$ is 
contained in the subspace 
$G_{\lambda, \mu}\cdot(\mathcal{S}_{\lambda}(x) \otimes \mathcal{S}_{\mu}(y))$, 
where the function $G_{\lambda, \mu}$ is defined by 
\bea 
G_{\lambda, \mu}&:=& 
\prod_{1 \le a<a' \le k}\prod_{i, j}
(x_{i}^{(a)}-x_{j}^{(a')})^{2a} 
\prod_{a=1}^{k}\prod_{i<j}
(x_{i}^{(a)}-x_{j}^{(a)})^{2a} \nn \\
&& {}\times 
\prod_{1 \le b<b' \le k}\prod_{i, j}
(y_{i}^{(b)}-y_{j}^{(b')})^{2b} 
\prod_{b=1}^{k}\prod_{i<j}
(y_{i}^{(b)}-y_{j}^{(b)})^{2b} \nn \\ 
&& {}\times 
\prod_{1 \le a, b \le k \atop a+b>k}\prod_{i, j}
(x_{i}^{(a)}-y_{j}^{(b)})^{a+b-k}
\prod_{a=1}^{k}\prod_{j}(x_{j}^{(a)})^{(a-b_{0})_{+}}. 
\ena 
\end{prop}

In the following, we prove that the image of 
$\varphi_{\lambda, \mu}|_{\Gamma_{\lambda, \mu}}$ is equal to 
$G_{\lambda, \mu}\cdot(\mathcal{S}_{\lambda}(x) \otimes \mathcal{S}_{\mu}(y))$ 
by using vertex operators in the same way as before.

\subsection{Construction of vertex operators} 

Decompose $\C^{2k}$ into $k$ orthogonal components 
\bea  
\C^{2k}=V_{1} \oplus \cdots \oplus V_{k}, 
\quad 
V_{j} \simeq \C^{2}, (j=1, \ldots , k). 
\ena 
We define a basis of $\C^{k}$ as follows. 
Take a basis $\{\epsilon_{j}^{+}, \epsilon_{j}^{-}\}$ 
of $V_{j}\simeq\C^{2}$ such that 
\bea 
\langle \epsilon_{j}^{\pm}, \epsilon_{j}^{\pm} \rangle=2, \quad 
\langle \epsilon_{j}^{\pm}, \epsilon_{j}^{\mp} \rangle=1. 
\label{defofepsilon}
\ena 
Then the set of vectors 
$\{\epsilon^{+}_{1}, \epsilon^{-}_{1}, \ldots , 
   \epsilon^{+}_{k}, \epsilon^{-}_{k}\}$ 
is a basis of $\C^{2k}$. 

Let 
$\boldsymbol\alpha_{j}=(\{\alpha_{j, m}\}, \alpha_{j}^{0})$ 
and 
$\boldsymbol\beta_{j}=(\{\beta_{j, m}\}, \beta_{j}^{0}), 
(j=1, \ldots , k)$ be 
sequences of vectors in $V_{j}\subset\C^{2k}$ defined by 
\bea 
\alpha_{j, m}=\alpha_{j}^{0}=\epsilon_{j}^{+}, \quad 
\beta_{j, m}=\beta_{j}^{0}=\epsilon_{j}^{-}, \quad 
(\forall m \in \Z). 
\ena 
We define the vertex operators 
$A_{a}(z)$ and $B_{b}(z), (a, b=1, \ldots , k)$ by 
\bea 
A_{a}(z):=X_{\boldsymbol\alpha_{a}}(z), \quad 
B_{b}(z):=X_{\boldsymbol\beta_{b}}(z). 
\ena 
These operators satisfy 
\bea 
&& 
A_{a}(z)A_{b}(w)=(z-w)^{2\delta_{a, b}}:\!A_{a}(z)A_{b}(w)\!:, \\ 
&&  
A_{a}(z)B_{b}(w)=(z-w)^{\delta_{a, b}}:\!A_{a}(z)B_{b}(w)\!:,  \\ 
&& 
B_{a}(z)A_{b}(w)=(z-w)^{\delta_{a, b}}:\!B_{a}(z)A_{b}(w)\!:, \\ 
&& 
B_{a}(z)B_{b}(w)=(z-w)^{2\delta_{a, b}}:\!B_{a}(z)B_{b}(w)\!:. 
\ena 
In particular, we have 
\bea 
A_{a}(z)A_{b}(w)=A_{b}(w)A_{a}(z), \quad 
B_{a}(z)B_{b}(w)=B_{b}(w)B_{a}(z) 
\label{commrel1:r=3'} 
\ena 
for $a, b=1, \ldots , k$, and 
\bea 
A_{a}(z)^{2}=0, \quad B_{b}(z)^{2}=0, \quad 
A_{a}(z)B_{a}(z)=0=B_{a}(z)A_{a}(z)
\label{commrel2:r=3'} 
\ena 
for $a=1, \ldots , k$. 

Now we set 
\bea 
&& 
A_{\varepsilon}(z):=\varepsilon_{1}A_{1}(z)+ \cdots +\varepsilon_{k}A_{k}(z), 
\label{defA:2k} \\ 
&& 
B_{\varepsilon}(z):=\varepsilon_{1}B_{k}(z)+ \cdots +\varepsilon_{k}B_{1}(z). 
\label{defB:2k} 
\ena 
Note that the ordering of operators is reversed in (\ref{defA:2k}) and (\ref{defB:2k}). 

Let $\lambda$ and $\mu$ be level-$k$ restricted partitions of $l_{1}$ and $l_{2}$, 
respectively. 
Define the vertex operatos $A_{\lambda}(x_{1}, \ldots , x_{n})$ 
and $B_{\mu}(y_{1}, \ldots , y_{l})$ by 
\bea 
&& 
A_{\lambda}(x_{1}, \ldots , x_{n}):=
\prod_{a=1}^{k}\frac{1}{\lambda_{a}'!} \Bigl.
\left(\frac{\partial}{\partial \varepsilon_{1}}\right)^{\lambda_{1}'} \cdots 
\left(\frac{\partial}{\partial \varepsilon_{k}}\right)^{\lambda_{k}'}
A_{\varepsilon}(x_{1}) \cdots A_{\varepsilon}(x_{n})
\Bigr|_{\forall \varepsilon_{a}=0}, \\
&& 
B_{\mu}(y_{1}, \ldots , y_{l}):=
\prod_{b=1}^{k}\frac{1}{\mu_{b}'!} \Bigl.
\left(\frac{\partial}{\partial \varepsilon_{1}}\right)^{\mu_{1}'} \cdots 
\left(\frac{\partial}{\partial \varepsilon_{k}}\right)^{\mu_{k}'}
B_{\varepsilon}(y_{1}) \cdots B_{\varepsilon}(y_{l})
\Bigr|_{\forall \varepsilon_{a}=0}, 
\ena 
where $\lambda'=(\lambda_{1}', \ldots , \lambda_{k}')$ 
and $\mu'=(\mu_{1}', \ldots , \mu_{k}')$ are the conjugates 
of $\lambda$ and $\mu$, respectively. 

Set 
\bea 
\epsilon_{\lambda, \mu}^{*}:= 
\sum_{a=1}^{k}\lambda'_{a}\epsilon_{a}^{+}+
\sum_{b=1}^{k}\mu'_{b}\epsilon_{k+1-b}^{-} \in \C^{2k}. 
\ena 
Let $\gamma_{0}$ be a vector in $\C^{2k}$ uniquely determined by 
\bea 
\langle \epsilon_{a}^{+}, \gamma_{0} \rangle = \left\{ 
\begin{array}{ll} 
0, & (a \le b_{0}), \\ 
1, & (a > b_{0}), 
\end{array}, \right. 
\quad 
\langle \epsilon_{b}^{-}, \gamma_{0} \rangle =0, \,\, (1 \le b \le k). 
\ena 
 Consider the space of functions 
\bea 
U_{\lambda, \mu}:=\{ 
\bra{\gamma_{0}+\epsilon_{\lambda, \mu}^{*}}
h A_{\lambda}(x_{1}, \ldots , x_{n}) B_{\mu}(y_{1}, \ldots , y_{l}) \ket{\gamma_{0}}; 
h \in \widehat{\mathcal{H}}_{2k}^{+} \}. 
\ena 
{}From (\ref{commrel1:r=3'}) it is easy to see that 
\bea 
U_{\lambda, \mu} \subset 
\C[x_{1}, \ldots , x_{n}]^{\mathfrak{S}_{n}} \otimes 
\C[y_{1}, \ldots , y_{l}]^{\mathfrak{S}_{l}}.  
\ena 
Moreover, in the same way as Proposition \ref{subset:r=2} we have 

\begin{prop} 
\bea 
U_{\lambda, \mu} \subset \Gamma_{\lambda, \mu}. 
\ena 
\end{prop} 

The image $\varphi_{\lambda, \mu}(U_{\lambda, \mu})$ is 
given as follows. 

\begin{prop} 
\bea 
\varphi_{\lambda, \mu}(U_{\lambda, \mu})= 
S_{\mathbf{m}(\lambda), \mathbf{m}(\mu)}
(\boldsymbol\gamma_{1}^{+}, \ldots , \boldsymbol\gamma_{k}^{+}, 
 \boldsymbol\gamma_{1}^{-}, \ldots , \boldsymbol\gamma_{k}^{-}; \gamma_{0}). 
\ena 
Here the right hand side is defined by (\ref{defofV}) 
with the substitution $x_{j}^{(k+b)}:=y_{j}^{(b)}, (b=1, \ldots , k)$. 
The sequences 
$\boldsymbol\gamma_{a}^{\pm}=(\{\gamma_{a, m}^{\pm}\}, \gamma_{a}^{\pm , 0})$
are given by 
\bea 
\gamma_{a, m}^{+}=\gamma_{a}^{+, 0}=\sum_{j=1}^{a}\epsilon_{j}^{+}, 
\quad 
\gamma_{a, m}^{-}=\gamma_{a}^{-, 0}=\sum_{j=1}^{a}\epsilon_{k+1-a}^{-}, 
\quad 
(\forall m \in \Z).
\ena  
\end{prop} 

Note that the vectors $\gamma_{a, -m}^{\pm}, (a=1, \ldots , k)$ 
are linearly independent for each $m>0$. 
Hence we can apply Theorem \ref{matrixelt}. 
The functions $g(z, w; \boldsymbol\gamma_{a}^{\pm}, \boldsymbol\gamma_{b}^{\pm})$ 
are given by 
\bea 
&& 
g(z, w; \boldsymbol\gamma_{a}^{\pm}, \boldsymbol\gamma_{b}^{\pm})
=(z-w)^{2{\rm min}(a, b)}, \\ 
&& 
g(z, w; \boldsymbol\gamma_{a}^{\pm}, \boldsymbol\gamma_{b}^{\mp})=\left\{ 
\begin{array}{lc} 
(z-w)^{a+b-k}, & {\rm if} \,\, a+b>k, \\ 
1, & {\rm if} \,\, a+b \le k,  
\end{array} \right. \\ 
&& 
\langle \gamma_{a}^{+, 0}, \gamma_{0} \rangle =(a-b_{0})_{+}, \quad 
\langle \gamma_{a}^{-, 0}, \gamma_{0} \rangle =0. 
\ena 
Therefore we see that 

\begin{cor}\label{result:r=3'}
\bea 
\varphi_{\lambda, \mu}(\Gamma_{\lambda, \mu})= 
G_{\lambda, \mu}\cdot(\mathcal{S}_{\lambda}(x) \otimes \mathcal{S}_{\mu}(y)). 
\ena 
\end{cor}

\subsection{Fermionic formula} 

{}From Proposition \ref{image:r=3'} and Corollary \ref{result:r=3'}, 
we have 
\bea 
\mathop{\rm ch}(\Gamma_{\lambda, \mu}/\Gamma_{\lambda, \mu}')= 
\mathop{\rm ch}(G_{\lambda, \mu}\cdot(\mathcal{S}_{\lambda}(x) \otimes \mathcal{S}_{\mu}(y))). 
\ena 

The character of 
$G_{\lambda, \mu}\cdot(\mathcal{S}_{\lambda}(x) \otimes \mathcal{S}_{\mu}(y))$ 
is given as follows. 
Introduce the $2k \times 2k$ matrix $A$ defined by 
\bea 
A:=\left( 
\begin{array}{c|c} 
A^{(2)} & B^{(3)} \\ 
\hline 
B^{(3)} & A^{(2)} 
\end{array} \right), 
\label{gordon-mat:another} 
\ena 
where $A^{(2)}$ is the matrix defined by (\ref{gordon-mat:r=2}) 
and $B^{(3)}$ is defined by 
\bea 
B^{(3)}=(B^{(3)}_{ab})_{1 \le a, b \le k}, \quad 
B^{(3)}_{ab}:={\rm max}(0, a+b-k). 
\label{gordon-mat:r=3} 
\ena 
For example, the matrix $A$ for $k=1, 2$ and $3$ is given by 
\bea 
\left( 
\begin{array}{cc} 
2 & 1 \\ 
1 & 2 
\end{array} 
\right), 
\left( 
\begin{array}{cccc} 
2 & 2 & 0 & 1 \\ 
2 & 4 & 1 & 2 \\
0 & 1 & 2 & 2 \\ 
1 & 2 & 2 & 4 
\end{array} 
\right) \,\, 
{\rm and} \,\, 
\left( 
\begin{array}{cccccc} 
2 & 2 & 2 & 0 & 0 & 1 \\ 
2 & 4 & 4 & 0 & 1 & 2 \\
2 & 4 & 6 & 1 & 2 & 3 \\ 
0 & 0 & 1 & 2 & 2 & 2 \\ 
0 & 1 & 2 & 2 & 4 & 4 \\ 
1 & 2 & 3 & 2 & 4 & 6 
\end{array} 
\right), 
\ena 
respectively. 
We denote by $\mathbf{c}^{(3)}_{b_{0}}$ the vector defined by 
\bea 
\mathbf{c}^{(3)}_{b_{0}}:=
(\underbrace{0, \ldots , 0, 1, 2, \ldots , k-b_{0}}_{k}, \underbrace{0, \ldots , 0}_{k}). 
\label{boundary:r=3} 
\ena 
Then we have 
\bea 
\mathop{\rm ch}(G_{\lambda, \mu}\cdot(\mathcal{S}_{\lambda}(x) \otimes \mathcal{S}_{\mu}(y)))= 
\frac{q^{\frac{1}{2}({}^{t}\mathbf{m} A \mathbf{m}-({\rm diag}A)\cdot \mathbf{m})
         +\mathbf{c}^{(3)}_{b_{0}}\cdot\mathbf{m}}}
     {(q)_{m_{1}(\lambda)} \cdots (q)_{m_{k}(\lambda)} (q)_{m_{1}(\mu)} \cdots (q)_{m_{k}(\mu)}}, 
\label{formula-GS:another} 
\ena 
where 
$\mathbf{m}:={}^{t}(m_{1}(\lambda), \ldots , m_{k}(\lambda), m_{1}(\mu), \ldots , m_{k}(\mu))$. 

{}From (\ref{formula-GS:another}) and Corollary \ref{fermionic:another}, 
we obtain the fermionic formula for the character of $(k, 3)$-admissible configurations: 
\begin{thm}\label{fermionic-formula:2k}
\bea 
&& 
\chi_{k, 3; (b_{0}, k)}(q, z) \\ 
&&{}=\sum_{n=0}^{\infty} \sum_{l_{1}+l_{2}=n \atop l_{1}, l_{2} \ge 0} 
\sum_{\sum_{j}jm_{i, j}=l_{i}, \atop (i=1, 2)}
\frac{q^{{}^{t}\mathbf{m} A \mathbf{m}-({\rm diag}A) \cdot \mathbf{m}
          +2\mathbf{c}_{b_{0}}^{(3)}\cdot \mathbf{m}}} 
     {(q^{2})_{m_{1,1}} \cdots (q^{2})_{m_{1,k}} (q^{2})_{m_{2, 1}} \cdots (q^{2})_{m_{2,k}}} 
q^{l_{2}} z^{n}, \nn
\ena 
where $A$ is the matrix defined by (\ref{gordon-mat:another}), 
$\mathbf{c}_{b_{0}}^{(3)}$ is the vector defined by (\ref{boundary:r=3}), 
$\mathbf{m}={}^{t}(m_{1, 1}, \ldots , m_{1, k}, m_{2, 1}, \ldots , m_{2, k})$ 
and $(q^{2})_{m}:=\prod_{j=1}^{m}(1-q^{2j})$. 
\end{thm}

\setcounter{section}{8}
\setcounter{equation}{0}

\section{Another fermionic formula for $\chi_{k, 3}$ in a special case} 

In this section we consider $(k, 3)$-admissible configurations 
with the initial condition $\mathbf{b}=([\frac{k+1}{2}], k)$, that is, 
$a_{0} \le [\frac{k+1}{2}]$. 
In this special case we can find another fermionic formula. 
As a consequence we get non-trivial equality 
between the different fermionic formulas for the character 
with $\mathbf{b}=([\frac{k+1}{2}], k)$. 

\subsection{Functional realization} 

First we give another functional realization of $E_{\Lambda}^{(k, 3)}$ 
for $\Lambda=b_{0}\Lambda_{1}+(k-b_{0})\Lambda_{2}$. 
We fix $\Lambda$ and abbreviate 
$E_{\Lambda}^{(k, 3)}$ and $J_{\Lambda}^{(k, 3)}$ 
to $E^{(k, 3)}$ and $J^{(k, 3)}$, respectively. 

Let $F_{n}=\C[x_{1}, \ldots , x_{n}]^{\mathfrak{S}_{n}}$. 
Define a pairing 
\bea 
\langle \cdot , \cdot \rangle: 
E_{n}^{(3)} \otimes F_{n} \longrightarrow \C 
\ena 
by 
\bea 
\langle e(z_{1}) \cdots e(z_{n}), f(x_{1}, \ldots , x_{n}) \rangle 
:=f(z_{1}, \ldots , z_{n}), 
\label{pairing:r=3res} 
\ena 
where 
\bea 
e(z):=e^{(1)}(z^{2})+ze^{(2)}(z^{2}). 
\ena 
This pairing is non-degenerate and respects the grading on $E_{n}^{(3)}$ 
and $F_{n}$. 

\begin{prop}
The orthogonal complement $F_{n}^{(k, 3)}:=(J_{n}^{(k, 3)})^{\perp}$ is 
the space of functions $f(x_{1}, \ldots , x_{n}) \in F_{n}$ such that 
\bea 
f=0 \quad {\rm if} \,\, 
\begin{array}{l} 
x_{1}=\cdots =x_{a}=-x_{a+1}=\cdots =-x_{k+1} \,\, (0\le \forall a \le k+1) \,\,{\rm or} \\ 
x_{1}=\cdots =x_{b_{0}+1}=0. 
\end{array} 
\ena 
\end{prop} 

\begin{proof} 
Recall that $J_{n}^{(k, 3)}$ is the ideal generated by the coefficients 
of $e^{(1)}(z)^{\alpha}e^{(2)}(z)^{\beta}, (\alpha+\beta=k+1)$ 
and the element $e^{(1)}(0)^{b_{0}+1}$. 
It is easy to see that the condition 
\bea 
e^{(1)}(z)^{\alpha}e^{(2)}(z)^{\beta}=0, \,\, {\rm for} \,\, \alpha+\beta=k+1 
\ena 
is equivalent to 
\bea 
e(z)^{a}e(-z)^{k+1-a}=0, \,\, {\rm for} \,\, 0 \le a \le k+1. 
\ena 
{}From this observation the proposition follows 
in the same way as Proposition \ref{perp-space}. 
\end{proof} 

Hence we have the following. 

\begin{prop} 
The pairing (\ref{pairing:r=3res}) induces a well-defined 
non-degenerate bilinear pairing of the graded spaces 
\bea 
\langle \cdot , \cdot \rangle: 
E_{n}^{(k, 3)} \otimes F_{n}^{(k, 3)} \longrightarrow \C, 
\ena 
where $E_{n}^{(k, 3)}$ is the graded component $E_{n}^{(k, 3)}:=E_{n}^{(3)}/J_{n}^{(k, 3)}$. 
\end{prop} 

Therefore the character of $(k, 3)$-admissible configurations is given as follows. 

\begin{cor} 
\bea 
\chi_{k, 3; (b_{0}, k)}=\sum_{n=0}^{\infty}z^{n}\mathop{\rm ch}F_{n}^{(k, 3)}(q). 
\ena 
\end{cor}

\subsection{Gordon filtration} 
For a level-$k$ restricted partition $\lambda$ 
we defined the map 
\bea 
\varphi_{\lambda}: 
\C[x_{1}, \ldots , x_{n}]^{\mathfrak{S}_{n}} \longrightarrow \mathcal{S}_{\lambda} 
\ena 
in (\ref{varphi:r=2}). 
Using this map we define the subspaces 
$\mathcal{F}_{\lambda}, \Gamma_{\lambda}$ and $\Gamma'_{\lambda}$
as in the case of $r=2$, 
that is, 
\bea 
\mathcal{F}_{\lambda} &:=& \mathop{\rm Ker} \varphi_{\lambda} \cap F_{n}^{(k, 3)}, \\ 
\Gamma_{\lambda} &:=& \cap_{\nu \succ \lambda} \mathcal{F}_{\nu}, \\
\Gamma'_{\lambda} &:=& \Gamma_{\lambda} \cap \mathop{\rm Ker} \varphi_{\lambda}.
\ena 
Then we have 
\bea 
\mathop{\rm ch}F_{n}^{(k, 3)}=\sum_{\lambda}\mathop{\rm ch}(\Gamma_{\lambda}/\Gamma_{\lambda}'). 
\label{decomp:r=3} 
\ena 

\begin{prop}\label{image:r=3res} 
Let $\lambda$ be a level-$k$ restricted partition of $n$. 
Then the image of the map $\varphi_{\lambda}|_{\Gamma_{\lambda}}$ is 
contained in the principal ideal $G_{\lambda}^{(3)}\mathcal{S}_{\lambda}$, 
where $G_{\lambda}^{(3)}$ is defined by 
\bea 
&& G_{\lambda}^{(3)}:=G_{\lambda}^{(2)}\overline{G}_{\lambda}^{(3)}, \\ 
&& \overline{G}_{\lambda}^{(3)}:= 
\prod_{1 \le a<b \le k \atop a+b>k}\prod_{i, j}
(x_{i}^{(a)}+x_{j}^{(b)})^{a+b-k} 
\prod_{a>\frac{k}{2}}\prod_{i<j}
(x_{i}^{(a)}+x_{j}^{(a)})^{2a-k}. 
\ena 
Here $G_{\lambda}^{(2)}$ is the function defined by (\ref{def-of-G:r=2}). 
\end{prop} 

\begin{proof} 
It suffices to prove that the function in the image of 
$\varphi_{\lambda}|_{\Gamma_{\lambda}}$ 
is divisible by $\overline{G}_{\lambda}^{(3)}$. 

Denote the variables $x_{p}$ such that $\varphi_{\lambda}(x_{p})=x_{i}^{(a)}$ 
by $x_{i, l}^{(a)}, (l=1, \ldots , a)$. 
We can carry out the evaluation of $\varphi_{\lambda}$ in two steps: 
$\varphi_{\lambda}(F)=\varphi_{2}(\varphi_{1}(F))$, 
where $\varphi_{1}$ is the evaluation of all variables except $\{x_{i, l}^{(a)}\}_{l=1}^{a}$ 
and $\varphi_{2}$ is the evaluation of the variables $\{x_{i, l}^{(a)}\}_{l=1}^{a}$. 
Let $F_{1}:=\varphi_{1}(F)$ for $F \in \Gamma_{\lambda}$. 
As a polynomial of $x_{i, l}^{(a)}, (l=1, \ldots , a)$, $F_{1}$ is symmetric. 
Moreover, $F_{1}$ equals zero if 
$(k-b+1)$ variables of $\{x_{i, l}^{(a)}\}$ are equal to $-x_{j}^{(b)}$ 
for $b=1, \ldots , k$ such that $a+b>k$. 
Therefore, the following lemma implies that 
$\varphi_{\lambda}(F)=\varphi_{2}(F_{1})$ is divisible 
by $\overline{G}_{\lambda}^{(3)}$. 
\end{proof} 

\begin{lem} 
Let $f(x_{1}, \ldots , x_{m})$ be a symmetric polynomial 
satisfying 
\bea 
f(x_{1}, \ldots , x_{m})=0, \quad {\rm if} \,\, x_{1}= \cdots =x_{s}=a 
\label{lem:asp:zerocond} 
\ena 
for some constant $a$. 
Then the polynomial $f(x, \ldots , x)$ is divisible by $(x-a)^{m-s+1}$. 
\end{lem} 

This lemma is easy to prove. 

%
%
%
%

If the induced map 
$\varphi_{\lambda}: \Gamma_{\lambda}/\Gamma_{\lambda}' \rightarrow 
                    G_{\lambda}^{(3)}\mathcal{S}_{\lambda}$
is surjective, we get the fermionic formula for $\mathop{\rm ch}F_{n}^{(k, 3)}$ 
by (\ref{decomp:r=3}). 
In fact the induced map is surjective. 
We can see this fact from the formula for 
the character of $(k, 3)$-admissible configurations obtained in \cite{combinatorial}. 
Here we do not assume this result. 
We have proved the surjectivity by using vertex operators 
only in the case of $b_{0}=[\frac{k+1}{2}]$. 
In the following we consider this special case.

\subsection{Construction of vertex operators}\label{VO:r=3}
 
First we consider the case that $k$ is odd. 
Set $k=2l+1$. 
Note that $b_{0}=[\frac{k+1}{2}]=l+1$. 

Decompose $\C^{k}$ into $l+1$ orthogonal components 
\bea 
\C^{k}=V \oplus V_{1} \oplus \cdots \oplus V_{l}, 
\quad 
V \simeq \C, \quad V_{j} \simeq \C^{2}, (j=1, \ldots , l). 
\ena 
Take a vector $\epsilon_{0} \in V\simeq\C$ such that 
$\langle \epsilon_{0}, \epsilon_{0} \rangle=1$. 
Next we take a basis $\{\epsilon_{j}^{+}, \epsilon_{j}^{-}\}$ 
of $V_{j}\simeq\C^{2}$ satisfying (\ref{defofepsilon}).
Then the set of vectors 
$\{\epsilon, \epsilon^{+}_{1}, \epsilon^{-}_{1}, \ldots , 
   \epsilon^{+}_{l}, \epsilon^{-}_{l}\}$ 
is a basis of $\C^{k}$. 

Let  
$\boldsymbol\alpha=(\{\alpha_{m}\}, \alpha^{0})$ 
be a sequence of vectors in $V \subset \C^{k}$ defined by 
\bea 
\alpha^{0}=\alpha_{2m}:=\sqrt{3}\epsilon_{0}, \quad 
\alpha_{2m+1}:=\epsilon_{0}, \quad 
(m \in \Z), 
\label{defofalpha'} 
\ena 
and  
$\boldsymbol\alpha^{\pm}_{j}=(\{\alpha^{\pm}_{j, m}\}, \alpha_{j}^{\pm, 0}), 
(j=1, \ldots , l)$ 
sequences of vectors in $V_{j}\subset\C^{k}$ defined by
\bea 
\alpha^{\pm}_{j, m}:=(\pm 1)^{m}\epsilon_{j}^{\pm}, \quad 
\alpha^{\pm, 0}_{j}:=\epsilon_{j}^{\pm}. 
\label{defofalpha} 
\ena 

We rename the sequences $\boldsymbol\alpha$ and 
$\boldsymbol\alpha_{j}^{\pm}, (j=1, \ldots ,l)$ 
to $\boldsymbol\beta_{a}, (a=1, \ldots , k)$ by 
\bea 
\boldsymbol\beta_{a}:=\left\{ 
\begin{array}{lc} 
\boldsymbol\alpha_{a}^{+}, & 1 \le a \le l, \\ 
\boldsymbol\alpha, & a=l+1, \\ 
\boldsymbol\alpha_{k-a+1}^{-}, & l+2 \le a \le k. 
\end{array} \right. 
\label{defofbeta} 
\ena 
Now we define the vertex operators $E_{a}(z), (a=1, \ldots , k)$ by 
\bea 
E_{a}(z):=X_{\boldsymbol\beta_{a}}(z). 
\label{defofE:k} 
\ena 
These operators satisfy the following: 
\bea 
E_{a}(z)E_{b}(w)=\left\{ 
\begin{array}{lc} 
(z-w)^{2}:E_{a}(z)E_{b}(w)\!:, & a=b \not= l+1, \\ 
(z-w)^{2}(z+w):E_{a}(z)E_{b}(w)\!:, & a=b=l+1, \\ 
(z+w):E_{a}(z)E_{b}(w)\!:, & a+b=k+1, a\not=l+1, \\ 
:E_{a}(z)E_{b}(w)\!:, & {\rm otherwise}. 
\end{array} \right.  
\ena 
In particular, 
\bea 
&& 
E_{a}(z)E_{b}(w)=E_{b}(w)E_{a}(z)
\label{commrel1}
\ena 
for $a, b=1, \ldots , k$, and 
\bea 
E_{a}(z)^2=0, \quad E_{a}(z)E_{k+1-a}(-z)=0 
\label{commrel2} 
\ena 
for $a=1, \ldots , k$. 
 
As in the case of $r=2$, 
we set 
\bea 
E_{\epsilon}(z):=\epsilon_{1}E_{1}(z)+ \cdots +\epsilon_{k}E_{k}(z). 
\ena 
For a level-$k$ restricted partition $\lambda$ of $n$, we set 
\bea 
E_{\lambda}(x_{1}, \ldots , x_{n}):=
\prod_{a=1}^{k}\frac{1}{\lambda_{a}'!} \Bigl.
\left(\frac{\partial}{\partial \epsilon_{1}}\right)^{\lambda_{1}'} \cdots 
\left(\frac{\partial}{\partial \epsilon_{k}}\right)^{\lambda_{k}'}
E_{t}(x_{1}) \cdots E_{t}(x_{n})
\Bigr|_{\forall \epsilon_{a}=0}, 
\label{defofElambda} 
\ena 
where $\lambda'=(\lambda_{1}', \ldots , \lambda_{k}')$ is the conjugate of $\lambda$. 

Set 
\bea 
\textstyle 
\epsilon_{\lambda}^{*}:=
\sum_{a=1}^{l}\lambda_{a}'\epsilon_{a}^{+}+\sqrt{3} \lambda_{l+1}'\epsilon_{0}+ 
\sum_{a=l+2}^{k}\lambda_{a}'\epsilon_{k+1-a}^{-} \in \C^{k}. 
\ena 
Then we see that 
\bea 
\bra{\beta} h E_{\lambda}(x_{1}, \cdots , x_{n}) \ket{0}=0, \, 
(\forall h \in \widehat{\mathcal{H}}_{k}^{+}), \quad 
{\rm unless} \,\, \beta=\epsilon_{\lambda}^{*}. 
\ena 
Consider the space of symmetric polynomials 
\bea 
U_{\lambda}:=
\{ \bra{\epsilon_{\lambda}^{*}} 
   h E_{\lambda}(x_{1}, \cdots , x_{n}) \ket{0}; 
   h \in \widehat{\mathcal{H}}_{k}^{+} \}. 
\label{defofUlambda} 
\ena 

\begin{prop}\label{inclusion:r=3}
\bea 
U_{\lambda} \subset \Gamma_{\lambda}.
\ena 
\end{prop} 

\begin{proof} 
In a similar way to the proof of Proposition \ref{subset:r=2}, 
it can be shown that $\varphi_{\mu}(U_{\lambda})=0$ for any $\mu$ such that $\mu > \lambda$. 
Hence it suffices to prove that $U_{\lambda} \subset F_{n}^{(k, 3)}$, 
and this is equivalent to 
\bea 
E_{\lambda}(\underbrace{x, \ldots , x}_{p}, \underbrace{-x, \ldots , -x}_{k+1-p}, 
            x_{k+2}, \ldots , x_{n})=0 
\ena 
for $p=0, \ldots , k+1$, and 
\bea 
E_{\lambda}(\underbrace{0, \ldots , 0}_{l+2}, x_{[(k+1)/2]+2}, \ldots , x_{n})=0. 
\ena 
This follows from the relation (\ref{commrel2}). 
\end{proof} 

{}From the relation $E_{a}(z)^{2}=0$, 
the following proposition holds as in the case of $r=2$: 

\begin{prop}\label{image:r=3}
\bea 
\varphi_{\lambda}(U_{\lambda})=
S_{\mathbf{m}(\lambda)}(\boldsymbol\gamma_{1}, \cdots , \boldsymbol\gamma_{k}). 
\ena 
Here the sequences of vectors 
$\boldsymbol\gamma_{a}=(\{\gamma_{a, m}\}, \gamma_{a}^{0}), (a=1, \ldots , k)$ 
are defined by 
\bea 
\textstyle
\gamma_{a, m}:=\sum_{j=1}^{a}\beta_{j, m}, \,\, (\forall m \in \Z), 
\quad \gamma_{a}^{0}:=\sum_{j=1}^{a}\beta_{j}^{0}, 
\ena 
where $\boldsymbol\beta_{j}=(\{\beta_{j, m}\}, \beta_{j}^{0}), (j=1, \ldots , k)$ 
are defined in (\ref{defofbeta}).  
\end{prop} 

Note that the vectors $\gamma_{a, -m}, (a=1, \ldots , k)$ are 
linearly independent for each $m>0$. 
Hence we can apply Theorem \ref{matrixelt}. 
Then the function $g(z, w; \boldsymbol\gamma_{a}, \boldsymbol\gamma_{b})$ 
is given by 
\bea 
g(z, w; \boldsymbol\gamma_{a}, \boldsymbol\gamma_{b})=\left\{ 
\begin{array}{lc} 
(z-w)^{2{\rm min}(a, b)}, & {\rm if} \,\, a+b \le k, \\ 
(z-w)^{2{\rm min}(a, b)}(z+w)^{a+b-k}, & {\rm if} \,\, a+b>k. 
\end{array} \right. 
\label{defofg} 
\ena 
Therefore we find 
\begin{cor}\label{result:r=3}
\bea 
\varphi_{\lambda}(\Gamma_{\lambda})= 
G_{\lambda}^{(3)}\mathcal{S}_{\lambda}.
\ena 
\end{cor}

Now we consider the case that $k$ is even. 
Set $k=2l$. 
The proof of Corollary \ref{result:r=3} for this case 
is quite similar to the case that $k$ is odd. 

First introduce the vertex operators $E_{a}(z), (a=1, \ldots , k)$ 
as follows. 

We decompose $\C^{k}$ into $l$ orthogonal components 
\bea  
\C^{k}=V_{1} \oplus \cdots \oplus V_{l}, 
\quad 
V_{j} \simeq \C^{2}, (j=1, \ldots , l). 
\ena 
Take a basis $\{\epsilon_{j}^{+}, \epsilon_{j}^{-}\}$ 
of $V_{j}\simeq\C^{2}$ satisfying (\ref{defofepsilon}).  
Let $\boldsymbol\alpha^{\pm}_{j}=(\{\alpha^{\pm}_{j, m}\}, \alpha_{j}^{\pm, 0}), 
(j=1, \ldots , l)$ be 
sequences of vectors in $V_{j}\subset\C^{k}$ defined by (\ref{defofalpha}). 
We rename the sequences 
$\boldsymbol\alpha_{j}^{\pm}, (j=1, \ldots ,l)$ 
to $\boldsymbol\beta_{a}, (a=1, \ldots , k)$ by 
\bea 
\boldsymbol\beta_{a}:=\left\{ 
\begin{array}{lc} 
\boldsymbol\alpha_{a}^{+}, & 1 \le a \le l, \\ 
\boldsymbol\alpha_{k-a+1}^{-}, & l+1 \le a \le k. 
\end{array} \right. 
\ena 
Then we define the vertex operators $E_{a}(z), (a=1, \ldots , k)$ by 
\bea 
E_{a}(z):=X_{\boldsymbol\beta_{a}}(z). 
\ena 
These operators satisfy the following: 
\bea 
E_{a}(z)E_{b}(w)=\left\{ 
\begin{array}{lc} 
(z-w)^{2}:E_{a}(z)E_{b}(w)\!:, & a=b, \\ 
(z+w):E_{a}(z)E_{b}(w)\!:, & a+b=k+1, \\ 
:E_{a}(z)E_{b}(w)\!:, & {\rm otherwise}. 
\end{array} \right.  
\ena 
The commutation relations (\ref{commrel1}) and (\ref{commrel2})
hold also in this case. 

Next we define the operator $E_{\lambda}(x_{1}, \ldots , x_{n})$ 
by (\ref{defofElambda}) for a level-$k$ restricted partition $\lambda$, 
and consider the space of matrix elements $U_{\lambda}$ 
defined by (\ref{defofUlambda}), 
where $\epsilon^{*}_{\lambda}$ is given by 
\bea 
\textstyle 
\epsilon^{*}_{\lambda}= 
\sum_{a=1}^{l}\lambda_{a}'\epsilon_{a}^{+}+ 
\sum_{a=l+1}^{k}\lambda_{a}'\epsilon_{k+1-a}^{-}. 
\ena 
Then it is easy to see that Proposition \ref{inclusion:r=3} 
and Proposition \ref{image:r=3} hold. 
The vectors $\gamma_{a, -m}, (a=1, \ldots , k)$ in Proposition \ref{image:r=3}  
are linearly independent for each $m>0$ also in this case. 
The function $g(z, w; \boldsymbol\gamma_{a}, \boldsymbol\gamma_{b})$ 
is given by (\ref{defofg}). 
Therefore Corolally \ref{result:r=3} holds also in the case that $k$ is even.

\subsection{Fermionic formula}

At last we write down the fermionic formula for 
$(k, 3)$-admissible configurations with the initial condition $a_{0} \le [\frac{k+1}{2}]$. 

{}From Proposition \ref{image:r=3res} and Colrollary \ref{result:r=3}, 
we have 
\bea 
\mathop{\rm ch}(\Gamma_{\lambda}/\Gamma_{\lambda}')=
\mathop{\rm ch}(G_{\lambda}^{(3)}\mathcal{S}_{\lambda}). 
\ena 
In order to write down the character of $G_{\lambda}^{(3)}\mathcal{S}_{\lambda}$ 
we introduce the $k \times k$ matrix $B$ defined by $B:=A^{(2)}+B^{(3)}$, 
that is, 
\bea 
B=(B_{ab})_{1 \le a, b \le k}, \quad B_{ab}:=2{\rm min}(a, b)+(a+b-k)_{+}. 
\label{Gordon:r=3res} 
\ena 
Then we have 
\bea 
\mathop{\rm ch}(G_{\lambda}^{(3)}\mathcal{S}_{\lambda})=
\frac{q^{\frac{1}{2}({}^{t}\mathbf{m} B \mathbf{m}
          -({\rm diag}B) \cdot \mathbf{m}) 
          +{\mathbf{c}}^{(2)}_{[\frac{k+1}{2}]}\cdot\mathbf{m} }} 
     {(q)_{m_{1}(\lambda)} \cdots (q)_{m_{k}(\lambda)}}, 
\ena 
where ${\mathbf{c}}^{(2)}_{[\frac{k+1}{2}]}$ is defined by (\ref{boundary-term:r=2}) 
with $b_{0}=[\frac{k+1}{2}]$. 

Finally we get 
\begin{thm} 
\bea 
\chi_{3, r; ([\frac{k+1}{2}], k)}(q, z)=\sum_{n=0}^{\infty} 
\sum_{m_{1}+2m_{2}+\cdots +km_{k}=n \atop m_{1}, \ldots , m_{k} \ge 0} 
\frac{q^{\frac{1}{2}({}^{t}\mathbf{m} B \mathbf{m}-({\rm diag}B) \cdot \mathbf{m}) 
           +\mathbf{c}^{(2)}_{[\frac{k+1}{2}]}\cdot\mathbf{m}}} 
     {(q)_{m_{1}} \cdots (q)_{m_{k}}} z^{n}, 
\ena 
where $B$ is the $k \times k$ matrix 
defined by (\ref{Gordon:r=3res}), 
$\mathbf{c}^{(2)}_{[\frac{k+1}{2}]}$ is the vector defined by (\ref{boundary-term:r=2}) 
and $\mathbf{m}={}^{t}(m_{1}, \ldots , m_{k})$. 
\end{thm}

\setcounter{section}{9} 
\setcounter{equation}{0} 

\section{Discussion} 

\subsection{}
The vertex operators constructed in Section \ref{VO:r=3} are a part of 
a vertex operator realization of $\widehat{\mathfrak{sl}}_{3}$ of level $k$. 
Here we describe the entire algebra $\widehat{\mathfrak{sl}}_{3}$ using the vertex operators 
in the cases of $k=1$ and $k=2$. 
For $k \ge 3$, the algebra is constructed as the tensor product 
of these algebras as mentioned in Introduction. 

{\it The $k=1$ case.} 
Set 
\bea 
&& 
\phi_{-}(z):=E_{1}(z)=X_{\boldsymbol\alpha}(z), \quad 
\phi_{+}(z):=X_{\boldsymbol{-\alpha}}(z), \\ 
&& 
\phi_{0}(z):=\sum_{n}a_{n}z^{-n-1}, \quad 
\overline{\phi}(z):= \, :\! \phi_{-}(-z)\phi_{+}(z) \!:, 
\ena 
where $\boldsymbol\alpha=(\{\alpha_{m}\}, \alpha^{0})$ is defined by (\ref{defofalpha'}) and 
$\boldsymbol{-\alpha}:=(\{-\alpha_{m}\}, -\alpha^{0})$. 
We abbreviated $a_{n}(\alpha_{n})$ to $a_{n}$. 
The operator product expansion is given as follows: 
\bea 
\phi_{-}(z)\phi_{+}(w) &\sim& \left\{ 
\begin{array}{l} 
\displaystyle 
\frac{1}{(z-w)^{2}}\frac{1}{2w}+\frac{1}{z-w}\left( \phi_{0}(w)-\frac{1}{2w} \right), 
\quad (z=w), \\ 
\displaystyle 
\frac{1}{z+w}\frac{\overline{\phi}(w)}{4w^{2}}, \quad (z=-w), 
\end{array} \right.  \\
\phi_{0}(z)\phi_{\pm}(w) &\sim& \left\{ 
\begin{array}{l}
\displaystyle  
\mp \frac{2\phi_{\pm}(w)}{z-w}, \quad (z=w), \\ 
\displaystyle 
\mp \frac{\phi_{\pm}(w)}{z+w}, \quad (z=-w), 
\end{array} \right. \\  
\phi_{0}(z)\overline{\phi}(w) &\sim& \left\{ 
\begin{array}{l} 
\displaystyle 
-\frac{\overline{\phi}(w)}{z-w}, \quad  (z=w), \\ 
\displaystyle 
\frac{\overline{\phi}(w)}{z+w}, \quad (z=-w), 
\end{array} \right. \\
\phi_{\pm}(z)\overline{\phi}(w) &\sim& \left\{ 
\begin{array}{l} 
\displaystyle 
\mp \frac{2w\phi_{\pm}(\pm w)}{z-w}, \quad (z=w), \\ 
\displaystyle 
0, \quad (z=-w) 
\end{array} \right. \\ 
\overline{\phi}(z)\overline{\phi}(w) &\sim& \left\{ 
\begin{array}{l} 
\displaystyle 
0, \quad (z=w), \\  
\displaystyle 
\frac{4w^{2}}{(z+w)^{2}}-\frac{4w}{z+w}\left( 1+2w(\phi_{0}(w)-\phi_{0}(-w)) \right), 
\,\, (z=-w), 
\end{array} \right. \\ 
\phi_{0}(z)\phi_{0}(w) &\sim& \left\{ 
\begin{array}{l} 
\displaystyle 
\frac{2}{(z-w)^{2}}, \quad (z=w), \\ 
\displaystyle 
-\frac{1}{(z+w)^{2}}, \quad (z=-w). 
\end{array} \right. 
\ena 
The operator $\phi_{\pm}(z)\phi_{\pm}(w)$ is regular at $z=\pm w$. 

The generators of $\widehat{\mathfrak{sl}}_{3}$ of level one 
are given by 
\bea 
&& 
\phi_{-}(z)=\sum_{n}e_{21}[n]z^{-2n}+\sum_{n}e_{31}[n]z^{-2n+1}, \\ 
&& 
\phi_{+}(z)=\sum_{n}e_{12}[n]z^{-2n-3}+\sum_{n}e_{13}[n]z^{-2n-4}, 
\ena 
\bea 
&&
-\frac{5}{4}+\frac{1}{2}z\phi_{0}(z)+\frac{1}{4}\overline{\phi}(z)= 
\sum_{n}e_{32}[n]z^{-2n+1}-\sum_{n}h_{13}[n]z^{-2n}, \\ 
&& 
-\frac{3}{4}+\frac{1}{2}z\phi_{0}(z)-\frac{1}{4}\overline{\phi}(z)= 
\sum_{n}e_{23}[n]z^{-2n-1}-\sum_{n}h_{12}[n]z^{-2n}.
\ena 
Here we set $h_{ij}:=e_{ii}-e_{jj}$. 

{\it The $k=2$ case.} 
Set 
\bea 
&& 
\phi_{-}(z):=E_{1}(z)+E_{2}(z)=
X_{\boldsymbol\alpha^{+}_{1}}(z)+X_{\boldsymbol\alpha^{-}_{1}}(z), \\  
&& 
\phi_{+}(z):=X_{\boldsymbol{-\alpha^{+}_{1}}}(z)+X_{\boldsymbol{-\alpha^{-}_{1}}}(z), \\ 
&& 
\phi_{0}(z):=\sum_{n}a_{n}z^{-n-1}, \\ 
&& 
\overline{\phi}(z):= \, :\!X_{\boldsymbol{\alpha^{+}_{1}}}(-z)X_{\boldsymbol{-\alpha^{-}_{1}}}(z)\!: 
                    +:\!X_{\boldsymbol{\alpha^{-}_{1}}}(-z)X_{\boldsymbol{-\alpha^{+}_{1}}}(z)\!:, 
\ena 
where $\boldsymbol\alpha^{\pm}_{1}$ is defined by (\ref{defofalpha}) 
and $a_{n}:=a_{n}(\alpha^{+}_{1, n})+a_{n}(\alpha^{-}_{1, n})$. 
Then the operator product expansion is given as follows: 
\bea  
\phi_{-}(z)\phi_{+}(w) &\sim& \left\{ 
\begin{array}{l} 
\displaystyle 
\frac{2}{(z-w)^{2}}+\frac{\phi_{0}(w)}{z-w}, 
\quad (z=w), \\ 
\displaystyle 
\frac{\overline{\phi}(w)}{z+w}, \quad (z=-w), 
\end{array} \right.  \\
\phi_{0}(z)\phi_{\pm}(w) &\sim& \left\{ 
\begin{array}{l}
\displaystyle  
\mp \frac{2\phi_{\pm}(w)}{z-w}, \quad (z=w), \\ 
\displaystyle 
\mp \frac{\phi_{\pm}(w)}{z+w}, \quad (z=-w), 
\end{array} \right. \\  
\phi_{0}(z)\overline{\phi}(w) &\sim& \left\{ 
\begin{array}{l} 
\displaystyle 
-\frac{\overline{\phi}(w)}{z-w}, \quad  (z=w), \\ 
\displaystyle 
\frac{\overline{\phi}(w)}{z+w}, \quad (z=-w), 
\end{array} \right. \\
\phi_{\pm}(z)\overline{\phi}(w) &\sim& \left\{ 
\begin{array}{l} 
\displaystyle 
\frac{\phi_{\pm}(\pm w)}{z-w}, \quad (z=w), \\ 
\displaystyle 
0, \quad (z=-w) 
\end{array} \right. \\ 
\overline{\phi}(z)\overline{\phi}(w) &\sim& \left\{ 
\begin{array}{l} 
\displaystyle 
0, \quad (z=w), \\  
\displaystyle 
-\frac{2}{(z+w)^{2}}+\frac{1}{z+w}\left( \phi_{0}(w)-\phi_{0}(-w) \right), 
\,\, (z=-w), 
\end{array} \right. \\ 
\phi_{0}(z)\phi_{0}(w) &\sim& \left\{ 
\begin{array}{l} 
\displaystyle 
\frac{4}{(z-w)^{2}}, \quad (z=w), \\ 
\displaystyle 
-\frac{2}{(z+w)^{2}}, \quad (z=-w). 
\end{array} \right. 
\ena 
The operator $\phi_{\pm}(z)\phi_{\pm}(w)$ is regular at $z=\pm w$. 

The generators of $\widehat{\mathfrak{sl}}_{3}$ of level two 
are given by 
\bea 
&& 
\phi_{-}(z)=\sum_{n}e_{21}[n]z^{-2n}+\sum_{n}e_{31}[n]z^{-2n+1}, \\ 
&& 
\phi_{+}(z)=\sum_{n}e_{12}[n]z^{-2n-2}+\sum_{n}e_{13}[n]z^{-2n-3}, \\
&&
-2+\frac{1}{2}\phi_{0}(z)+\frac{1}{2}\overline{\phi}(z)= 
\sum_{n}e_{32}[n]z^{-2n}-\sum_{n}h_{13}[n]z^{-2n-1}, \\ 
&& 
-1+\frac{1}{2}\phi_{0}(z)-\frac{1}{2}\overline{\phi}(z)= 
\sum_{n}e_{23}[n]z^{-2n-2}-\sum_{n}h_{12}[n]z^{-2n-1}.
\ena 

\subsection{}
Our problem is to obtain the fermionic formula for 
the character of $(k, r)$-admissible configurations 
with the initial condition (\ref{boundary-admissible}). 
In previous sections we obtained the fermionic formulas 
for $(k, 2)$ and $(k, 3)$-admissible cofigurations with the condition $a_{0} \le b_{0}$. 
For the case of $r=2$ our result is sufficient because the condition $a_{0} \le b_{0}$ 
is the only initial condition. 
However, in the case of $r=3$, we should consider not only the condition $a_{0} \le b_{0}$ 
but $a_{0}+a_{1} \le b_{1}$. 
The fermionic formula we obtained in Section \ref{section:2k} 
is for the case of $b_{1}=k$. 
Here we consider the case of $b_{1} <k$. 

The definition (\ref{realization:2k}) of the space $I_{l_{1}, l_{2}}^{(k, 3)}$ 
is replaced by 
\bea 
g=0 \quad {\rm if} && 
x_{1}= \cdots =x_{a}=y_{1}= \cdots =y_{b}, \,\,  
(a \ge 0, b \ge 0, a+b=k+1), \nn \\ 
&& {\rm or} \,\,\, 
x_{1}= \cdots =x_{b_{0}+1}=0, \\ 
&& {\rm or} \,\,\, 
x_{1}= \cdots =x_{s}=y_{1}= \cdots =y_{t}=0, \,\, 
(s \ge 0, t \ge 0, s+t=b_{1}+1). \nn
\ena 
The functional realization $F_{l_{1}, l_{2}}^{(k, 3)}$ is given by 
(\ref{defofF:2k}) with this redefined space $I_{l_{1}, l_{2}}^{(k, 3)}$. 

Now introduce the filtration $\{ \Gamma_{\lambda, \mu} \}$ on $I_{l_{1}, l_{2}}^{(k, 3)}$ 
by (\ref{gordon-filt:2k}) 
and consider the image of $\varphi_{\lambda, \mu}|_{\Gamma_{\lambda, \mu}}$ 
as in Proposition \ref{image:r=3'}. 
Then the image is contained a space of functions described as follows. 
For a partition $\rho=(\rho_{1}, \rho_{2}, \ldots )$ 
denote by $m_{\rho}(x_{1}, \ldots , x_{n})$ 
the monomial symmetric funtion: 
\bea 
m_{\rho}(x_{1}, \ldots , x_{n}):={\rm Sym}(x_{1}^{\rho_{1}} \cdots x_{n}^{\rho_{n}}). 
\ena 
Let $I_{\lambda, \mu}$ be the ideal of $\mathcal{S}_{\lambda}(x)\otimes \mathcal{S}_{\mu}(y)$ 
generated by the elements 
\bea 
m_{\rho^{(1)}}(x_{1}^{(a)}, \ldots , x_{m_{a}(\lambda)}^{(a)}) 
m_{\rho^{(2)}}(y_{1}^{(b)}, \ldots , y_{m_{b}(\mu)}^{(b)}) 
\ena 
such that 
\bea 
b_{1}<a+b \le k, \quad m_{a}(\lambda)\not= 0,  \quad m_{b}(\mu) \not=0 
\ena 
and 
\bea  
\rho_{m_{a}(\lambda)}^{(1)}+\rho_{m_{b}(\lambda)}^{(2)} \ge {\rm min}(a, b-(b_{1}-b_{0})). 
\ena 
Then it can be shown that 
\bea 
\varphi_{\lambda, \mu}(\Gamma_{\lambda, \mu}) \subset 
G_{\lambda, \mu} \prod_{b=1}^{k}\prod_{j}(y_{j}^{(b)})^{(b-b_{1})_{+}} \cdot I_{\lambda, \mu}. 
\label{image:r=3:gengeral} 
\ena 
In Section \ref{section:2k} we proved that the two spaces in (\ref{image:r=3:gengeral}) 
are equal in the case of $b_{1}=k$ using the vertex operators. 
For the case of $b_{1}<k$ we do not have proof or disproof of this equality. 
\newline

\noindent
{\it Acknowledgments.}\quad 
BF is partially supported by the grants 
RFBR 02-01-01015 and INTAS-00-00055. 
JM is partially supported by 
the Grant-in-Aid for Scientific Research (B2) no.12440039, 
and TM is partially supported by 
(A1) no.13304010, Japan Society for the Promotion of Science.
The work of EM is partially supported by NSF grant DMS-0140460. 
YT is supported by the Japan Society for the Promotion of Science.


\begin{thebibliography}{99}

\bibitem[FJLMM1]{bosonic1} 
B. Feigin, M. Jimbo, S.Loktev, T. Miwa and E. Mukhin, 
{\it Bosonic formulas for $(k, l)$-admissible partitions}, 
math.QA/0107054. 

\bibitem[FJLMM2]{bosonic2} 
B. Feigin, M. Jimbo, S.Loktev, T. Miwa and E. Mukhin, 
{\it Addendum to 'Bosonic formulas for $(k, l)$-admissible partitions'}, 
math.QA/0112104. 

\bibitem[FJMM1]{jack} 
B. Feigin, M. Jimbo, T. Miwa and E. Mukhin, 
{\it A differential ideal of Symmetric polynomials spanned 
by Jack polynomials at $\beta=-(r-1)/(k+1)$}, 
Int. Math. Res. Notices, {\bf 23} (2002), 1223-1237. 

\bibitem[FJMM2]{macdonald} 
B. Feigin, M. Jimbo, T. Miwa and E. Mukhin, 
{\it Symmetric polynomials vanishing on the shifted diagonals 
and Macdonald polynomials}, 
math.QA/0209042. 

\bibitem[FJMMT1]{root-of-unity} 
B. Feigin, M. Jimbo, T. Miwa, E. Mukhin and Y. Takeyama, 
{\it Symmetric polynomials vanishing on the diagonals shifted by roots of unity}, 
math.QA/0209126.  

\bibitem[FJMMT2]{combinatorial} 
B. Feigin, M. Jimbo, T. Miwa, E. Mukhin and Y. Takeyama, 
{\it Particle content of the $(k, 3)$-configurations}, 
to appear.  

\bibitem[FK]{FK}
I. B. Frenkel and V. G. Kac, 
{\it Basic representations of affine Lie algebras and dual resonance models}, 
Invent. Math., {\bf 62} (1980), 23-66. 

\bibitem[FKLMM]{coinv3} 
B. Feigin, R. Kedem, S.Loktev, T. Miwa and E. Mukhin, 
{\it Combinatorics of the $\widehat{sl}_{2}$ spaces 
of coinvariants III}, 
math.QA/0012190. 

\bibitem[FS]{FS} B. Feigin and A. Stoyanovsky, {\it Quasi-particles
models for the representations of Lie algebras and geometry of flag
manifold}, hep-th/9308079, RIMS 942; {\it Functional models for the
representations of current algebras and the semi-infinite Schubert
cells}, Funct. Anal. Appl. {\bf 28} (1994), 55--72.

\bibitem[P]{P} 
M. Primc, 
{\it Vertex operator construction of standard modules for $A_{n}^{(1)}$}, 
Pacific J. Math. {\bf 162} (1994), 143-187. 

\end{thebibliography}

\end{document}